\documentclass[11pt,twoside,reqno]{amsart}

\usepackage{amssymb}
\usepackage{xcolor}
\usepackage{hyperref}
\usepackage[dvipsnames]{xcolor}
\usepackage{enumitem}
\usepackage{amsmath}
\usepackage{theoremref}
\usepackage{mathtools}
\usepackage{dsfont}
\usepackage{float}

\hypersetup{unicode=false, colorlinks=true, linkcolor=RoyalBlue, anchorcolor=RoyalBlue, citecolor=ForestGreen, filecolor=red, menucolor=RoyalBlue, urlcolor=RoyalBlue}

\def\const{\text{\rm const}}

\def\ti{\tilde}

\def\PW{\text{\rm PW}}
\def\dist{\text{\rm dist}}
\def\supp{\text{\rm supp}\,}
\def\to{\rightarrow}

\def\si{{\rm sinc}}
\def\sign{\text{\rm sign}}

\def\no{\noindent}

\def\R{{\mathbb R}}
\def\T{{\mathbb T}}

\def\Z{{\mathbb{Z}}}
\def\C{{\mathbb{C}}}
\def\N{{\mathbb{N}}}

\def\PP{{\mathcal{P}}}

\def\BB{{\mathcal B}}
\def\MM{{\mathcal M}}

\def\FF{{\mathcal F}}
\def\SS{{\mathcal S}}
\def\HH{{\mathcal H}}
\def\KK{{\mathcal K}}

\def\EE{{\mathcal E}}

\def\e{\varepsilon}
\def\d{\delta}

\def\L{\Lambda}
\def\l{\lambda}

\def\a{\alpha}
\def\b{\beta}
\def\s{\sigma}
\def\lan{\lambda_n}

\def\const{\text{\rm const}}

\def\sign{\text{\rm sign}}

\def\dist{\text{\rm dist}}
\def\supp{\text{\rm supp}\,}

\def\Sch{Schr\"odinger }
\def\ONP{\text{\rm ONP}}

\DeclareMathOperator{\diagop}{diag}
\newcommand{\diag}[1]{\diagop\!\left(#1\right)}

\theoremstyle{plain}
\newtheorem{theorem}{Theorem}[section]
\newtheorem{lemma}[theorem]{Lemma}
\newtheorem{corollary}[theorem]{Corollary}
\newtheorem{proposition}[theorem]{Proposition}
\newtheorem{remark}[theorem]{Remark}

\newtheorem{example}[theorem]{Example}
\newtheorem{counterexample}[theorem]{Counterexample}
\newtheorem{prop}[theorem]{Proposition}

\numberwithin{equation}{section}

\title{Problems in spectral analysis of canonical Hamiltonian systems}

\begin{document}

\dedicatory{
        To Barry Simon, our colleague and teacher, on the occasion of his 80's birthday.}
	
\author{N.~Makarov}
\address{California Institute of Technology\\
    Department of Mathematics\\
    Pasadena, CA 91125, USA}
\email{makarov@its.caltech.edu}

\author{A.~Poltoratski}
\address{University of Wisconsin\\ Department of Mathematics\\ Van Vleck Hall\\
    480 Lincoln Drive\\
    Madison, WI  53706\\ USA }
\email{poltoratski@wisc.edu}
\thanks{The second author was partially supported by
    NSF Grant DMS-2244801.}
    
\author{A.~Zhang}
\address{Vanderbilt University\\ Department of Mathematics\\ 1326 Stevenson Center\\
    Nashville, TN  37240\\ USA }
\email{ashley.zhang@vanderbilt.edu}

\begin{abstract} 
This note focuses on recent results in spectral analysis of canonical systems of differential equations obtained via the approach developed in our previous papers \cite{MIF1, MP3, etudes, etudes2, PZ, Direct}. 
Many of our results are motivated by the pioneering research of Barry Simon and his co-authors; 
see, for instance, the papers cited in the main text.
We discuss solutions to the inverse spectral problem (ISP) for canonical Hamiltonian systems and mixed spectral problems for Schr\"odinger operators. One of our goals is to show connections of ISP with classical tools of analysis, such as the Hilbert transform, orthogonal polynomials, the gap problem and solutions to the Riemann-Hilbert problem. We illustrate our results with  examples and discuss further questions.
\end{abstract}

\maketitle

\section{Introduction}

This note focuses on various aspects of spectral problems for second order differential operators on the half-line.

A {\it regular} half-line canonical (Hamiltonian) system is the equation
\begin{equation}
	\Omega\dot X = z \HH X \qquad {\rm on} \quad [0,\infty).\label{eq001A}
\end{equation}
Here the Hamiltonian $\HH = \HH(t)$ is a given $2\times 2$ matrix-function satisfying $$\HH \in  L^1_{\rm loc}[0,\infty), \quad \HH \ne 0 \quad {\rm a.e.}, \quad \HH \ge 0 \quad {\rm a.e.}$$
The first relation means that the entries of $\HH$ are integrable on each finite interval. Systems satisfying this condition are called {\it regular}. 
The matrix $\Omega$ in \eqref{eq001A} is the symplectic matrix $$\Omega=\begin{pmatrix} 0 & 1 \\ -1 & 0 \end{pmatrix},$$ and $z \in \C$ is the 'spectral parameter'.
The unknown function $X = X(t,z)$ is a two-dimensional vector-function on $[0,\infty)$.

The Krein-de Branges theory translates spectral problems for canonical systems in the language of complex analysis. It can be applied to a wide range of problems in mathematical physics, including Schr\"odinger operators, Dirac systems, and string equations. For a recent account of the theory see \cite{Rem}.



Classical spectral problems are divided into two main groups: direct and inverse problems. In direct problems one needs to find the spectra or the spectral measure of the operator from its Hamiltonian $\HH$ (or potential $q$), while in inverse problems one aims to recover $\HH$ (or $q$) from the spectral information.
A famous result by Marchenko \cite{M1, M2} says that $q$ can be uniquely recovered from the spectral measure. 
Another classical result is Borg's two-spectra theorem \cite{Borg}, which says that $q$ can be recovered from two spectra, corresponding to different pairs of boundary conditions, see Section \ref{functions} for further discussion.

In the first part of this paper (Sections~\ref{prelim}-\ref{Hom}), we focus on spectral problems for canonical Hamiltonian systems on the half-line, with an emphasis on the inverse problem.
We extend the classical results of Marchenko, Gelfand-Levitan, and Krein to broader classes of canonical systems, providing explicit solution algorithms, convergence results, and formulas, as well as explicit examples. 
This framework also allows us to investigate connections between inverse spectral problems and classical tools of analysis, including the Hilbert transform, orthogonal polynomials, the solution of Riemann–Hilbert problems and Bessel functions. 
Throughout our analysis, we emphasize computational aspects, providing explicit formulas for Hamiltonian recovery.

A relatively new type of spectral problems, the so-called mixed spectral problem, asks to recover the operator from partial information on the potential and the spectrum.
The well-known theorem by Hochstadt and Lieberman from 1978 \cite{HL} says that for a Schr\"odinger equation on $(0,\pi)$ knowing the potential on one-half of the interval $(0,\pi/2)$
and knowing one of the spectra allows one to recover the operator uniquely. 
The result is precise in the sense that the knowledge of the spectrum minus one point, or of the potential on $(0,\pi/2-\e)$, is insufficient. Further results in the same direction obtained by del Rio, Gesztesy, Horvath, Simon and others \cite{S, SG1, DGS, Horvath, MIF1} allow one to combine various parts of the spectral and direct information to recover the operator.

Motivated by Simon and Gesztesy's counterexample in \cite{SG2}, the second part of the paper (Sections \ref{PrelimUncertainty}–\ref{uncertainty}) combines ideas from the theory of mixed spectral problems for differential operators with recent developments in the Uncertainty Principle (UP) in harmonic analysis.
Using new results on Gap and Type problems of UP, we present a version of Borg’s two-spectra theorem for Schrödinger operators that allows for uncertainty in the locations of the eigenvalues. 
We also provide an explicit formula for the exact “size of uncertainty,” expressed in terms of the lengths of the intervals in which the eigenvalues are permitted to lie. 
Among other applications, we describe pairs of indeterminate operators arising in the three-interval case of the mixed spectral problem.

The paper is organized as follows:

Section~\ref{prelim} surveys the basics of the Krein-de Branges theory and the necessary results and formulas used in spectral problems.
Section~\ref{Etudes} focuses on solving the inverse spectral problem for Paley-Wiener systems.
Section~\ref{periodic} presents an explicit algorithm for periodic spectral measures, as well as convergence results for non-periodic measures.
Section~\ref{Hom} studies canonical systems with homogeneous and quasi-homogeneous spectral measures.
Section~\ref{PrelimUncertainty} contains additional background material for mixed spectral problems.
Section~\ref{GapType} discusses gap and type problems, and the connection between mixed spectral problems for \Sch operators and the Beurling--Malliavin problem on completeness of exponential systems.
Section~\ref{uncertainty} presents a version of the two-spectra theorem with uncertainty and discusses the three-interval case.
Section~\ref{examples} concludes with a variety of examples illustrating our methods and formulas, as well as several open problems.

Throughout the paper, we emphasize ideas and structural results; in most places, we omit detailed calculations and complete proofs, referring the readers to the cited literature for full arguments.


\section{Preliminaries to spectral problems}\label{prelim}

\subsection{Canonical systems and de Branges spaces}\label{CSdB} Instead of a two-dimensional vector function $X$ discussed in the introduction, one may look for a $2 \times 2$ matrix-valued solution $M = M(t,z)$ solving \eqref{eq001A}.
Such a matrix valued function satisfying the initial condition $M(0,z)=I$ is called the \textit{transfer matrix} or the \textit{matrizant} of the system. The columns of the transfer matrix $M$ are the solutions for the system \eqref{eq001A} satisfying the initial conditions $\begin{pmatrix} 1 \\ 0 \end{pmatrix}$ (Neumann) and $\begin{pmatrix} 0 \\ 1 \end{pmatrix}$ (Dirichlet) at $0$.
As a general rule we denote
\begin{equation} M = \begin{pmatrix} A & B \\ C & D \end{pmatrix}.\label{eqTM} \end{equation}

An entire function $F(z)$ belongs to the Hermite-Biehler (HB) class if 
$$|F(z)| > |F(\overline z)|\textrm{ for all } z \in \C_+.$$
We say that an entire function is real if it is real on $\R$.

For each fixed $t$, the entries of the transfer matrix $M$ of the system \eqref{eq001A}, $A(z) = A(t,z) \equiv A_t(z)$, $B(z), C(z)$ and $D(z)$ are real entire functions.
The functions $$E:= A - iC, \qquad \tilde E:= B - iD$$ belong to the Hermite-Biehler class defined above; see for instance \cite{Rem}.

For an entire function $G$ we denote by $G^\#$ its Schwarz reflection with respect to $\R$,
$G^\#(z) = \overline G(\overline z)$. We denote by $H^2$ the standard Hardy space in the upper half-plane.

For every Hermite-Biehler (HB) function $F$ one can consider the de Branges (dB) space $\BB(F)$, a Hilbert space of entire functions defined as
$$ \BB(F) = \left\{ G \ | G \text{ is entire,} \frac GF, \ \frac {G^\#}F \in H^2 \right\}. $$
The Hilbert space structure in $\BB(F)$ is inherited from $H^2$:
$$ <G, H>_{\BB(F)} = \left< \frac GF,\frac HF \right>_{H^2} = \int_{-\infty}^{\infty} G(t)\overline H(t) \frac{dt}{|F(t)|^2}.$$

The space $\BB(E)$ is a reproducing kernel Hilbert space, i.e., for each $\l \in \C$ there exists $K(\l, \cdot) \in \BB(E)$ such that for any $F \in \BB(E)$,
$$F(\l) = <F, K(\l, \cdot)>_{\BB(E)}.$$
The function $K(\l,z)$ is called the reproducing kernel (reprokernel) for the point $\l$.
In the case of the dB-space $\BB(E)$, $K(\l,z)$ has the formula
$$ K(\l,z) = \frac{1}{2\pi i} \frac{ E(z) E^\#(\overline \l)- E^\#(z)E(\overline \l)}{\overline \l -z } = \frac{1}{\pi} \frac{A(z) C(\overline \l) - C(z)A(\overline \l))}{\overline \l-z}, $$
where $A = (E + E^\#)/2$ and $C = (E^\# - E)/2i$ are real entire functions such that $E = A-iC$.

The functions $E,\ \ti E$ corresponding to a canonical system \eqref{eq001A} give rise to the family of dB-spaces
$$ \BB_t = \BB( E(t,\cdot) ), \qquad \tilde \BB_t = \BB( \tilde E( t,\cdot) ).$$

A value $t$ is $\HH$-regular if it does not belong to an open interval on which $\HH$ is a constant matrix of rank one.
The spaces $\BB_t, \ti \BB_t$ form {\it chains}, i.e., $\BB_s\subsetneq \BB_t$ for $s < t$ and the inclusion is isometric for regular $t$ and $s$.

If $E$ is an HB-function non-vanishing on $\R$, then one
can define a continuous branch of the argument of 
$E$. The function $\phi_E(x)=-\arg E(x)$ is called the phase function for the space $\BB(E)$. It is not difficult to show that $\phi_E$ is a growing
real analytic function, $\phi'>0$ on $\R$.

Special role in our formulas is played by the kernels at 0:  $K_t(z,0)$ and $\tilde K_t(z,0)$.
We denote $k_t(z) = K_t(z,0)$.

%
%

\subsection{Spectral measures}\label{secSM} 
There are several ways to introduce spectral measures of canonical systems. 
We'll make a simplifying assumption that the system has no "jump intervals", i.e., intervals on which the Hamiltonian is rank one and constant. 
In this case all $t\in [0,\infty)$ are $\HH$-regular and all inclusions $\BB_s\subset\BB_t, \ti\BB_s\subset\ti\BB_t$ are isometric.
We can make this assumption because we will be mostly concerned with the case $\det\HH\neq 0$ a.e.


A measure $\mu$ on $\R$ is called Poisson-finite ($\Pi$-finite) if  $$\int \frac{d|\mu|(x)}{1 + x^2} < \infty.$$
A measure on $\hat{\mathbb{R}}=\R\cup \{\infty\}$ is $\Pi$-finite if it is the sum of a $\Pi$-finite measure on $\R$ and a finite point mass at infinity.

For a dB-space $\BB(E)$ there exists a standard one-parameter family of discrete measures on $\R$ 
such that $\BB(E)\overset{\rm iso}= L^2(\mu)$. 
Let
$\a\in \C,\ |\a|=1$ and denote by
$t_n$ the sequence of points on $\R$ such that
$$\frac{E^\#(t_n)}{E(t_n)}=\alpha.$$
In terms of the phase function, $\phi(t_n)=\arg\alpha\ (\rm{mod }\ \pi)$.
Consider the measure
\begin{equation}\label{eqRM}
\mu_\alpha=\sum_n \frac \pi{\phi'(t_n)|E(t_n)|^2}\d_{t_n}.
\end{equation}
Then $\BB(E)\overset{\rm iso}= L^2(\mu_\alpha)$, for all
values of $\alpha$ except possibly one, see for instance \cite{dB}.

By definition, a positive measure $\mu$ on $\R$ is a spectral measure of the CS \eqref{eq001A} with the initial condition $\begin{pmatrix} 1 \\ 0 \end{pmatrix}$ at $t = 0$ if
$$\forall t,\qquad \BB_t\stackrel{\rm iso}{\subset} L^2(\mu).$$
(The definition is slightly more complicated in presence of jump intervals for the Hamiltonian, which we do not allow in this paper.)
It is well-known that spectral measures of regular CS are $\Pi$-finite; see for instance \cite{Rem}. 
In a similar way, using $\ti B_t$, one can define a spectral measure
$\ti \mu$ for the initial condition  $\begin{pmatrix}0\\1\end{pmatrix}$. 

Conversely, one of the main results of the Krein-de Branges theory says that every positive $\Pi$-finite measure is a spectral measure of a regular CS.
Note that in general, the corresponding system may not satisfy $\det \HH\neq 0$ a.e., the restriction we are assuming in the first part of this article.
Also, HB functions corresponding to the systems considered in this paper have no zeros on the real line.
We will assume this restriction in our general discussions of dB-spaces. 
(If $E$ vanishes at some point of $\R$, then all functions in $\BB(E)$ must vanish at the same point, as follows from the definition. PW type spaces discussed in this note clearly do not have such a property.)

Every regular canonical system has a spectral measure; in fact for $\mu$ we can take any limit point of the family of measures $|E_t|^{-2}$ as $t\to \infty$, \cite{dB}. The spectral measure may or may not be unique. 
It is unique iff
\begin{equation*}\int{\rm trace}\ \HH(t)dt=\infty.\end{equation*}
The case when the spectral measure is unique is called the limit point case and the case when it is not, the limit circle case.

Finally, let us mention that spectral measures are invariant with respect to 'time' parameterizations, i.e., a change of variable $t$ in the initial system \eqref{eq001A} via an increasing homeomorphism $t \mapsto s(t)$ does not
change the spectral measure.

\subsection{Paley-Wiener spaces}\label{secPW}
We will use the following definition for the Fourier transform in $L^2(\R)$:
$$ (\FF f)(\xi) \equiv \hat f(\xi) = \frac1{\sqrt{2\pi}} \int e^{-i\xi x}f(x)dx,$$
first defined on test functions and then extended to a unitary operator
$\FF L^2(\R)\to L^2(\R)$ via Parseval's theorem.
The Paley-Wiener space $\PW_t$ of entire functions is defined as the image
$${\rm PW}_t = \FF L^2[-t,t].$$
By the Paley-Wiener theorem, $\PW_t$ can be equivalently defined as the space of entire functions of exponential type at most $t$ which belong to $L^2(\R)$.
The Hilbert structure in $\PW_t$ is inherited from $L^2(\R)$.

Paley-Wiener spaces appear as the de Branges chain $\BB_t$ for the free system $\HH=I$.

\subsection{Aleksandrov-Clark and dual measures}\label{secClark}
Let $\varphi$ be a function from the unit ball of $H^\infty(\C_+)$
(such functions are often called Schur functions).
Associated with every such $\varphi$ is the family of $\Pi$-finite positive  measures $\{\s^\varphi_\a\}_{\a\in\T}$ on $\hat\R$ defined via the formula
\begin{equation} \label{eqClark}
\Re \frac{\a+\varphi(z)}{\a-\varphi(z)} = 
py + \frac{1}{\pi}\int{\frac{yd\s^\varphi_\a(t)}{(x-t)^2+y^2}},
\hspace{1cm} z=x+iy,
\end{equation}
where the number $\pi p $ can be interpreted as the pointmass of $\s^\varphi_\a$ at $\infty$.
Such measures $\s^\varphi_a$ are called Aleksandrov-Clark (AC) measures for $\varphi$.
Families of AC measures appear in complex function theory and applications to perturbation theory and spectral problems, see for instance \cite{PS}.

At most one of the measures $\{\sigma^\varphi_\alpha\}_{\alpha\in\T}$ has a pointmass at infinity:
\[
\sigma_\alpha(\infty)\ne 0 \text{ iff } \varphi(\infty)=\alpha,
\]
in the sense of non-tangential limit, and $\varphi-\alpha\in L^2(\R)$.

It follows from the definition that all of the measures $\s^\varphi_\a$ are singular
iff one of them is singular
iff $\varphi$ is inner.
All measures are discrete 
iff one of them is discrete 
iff $\varphi$ is a meromorphic inner function (MIF), an inner function which can be extended meromorphically to the whole plane.

If $\varphi=\theta$ is a MIF, then
\begin{equation}\label{eqClarkpp}
    \sigma_\alpha\cdot 1_\R=2\pi\sum_{\theta(\xi)=\alpha}\frac{\delta_\xi}{|\theta'(\xi)|},
\end{equation}
as can be deduced from the definition.

Every HB function $E(z)$ gives rise to a MIF, $\theta_E = E^\#/E$.
In this case the last formula is closely
related to \eqref{eqRM} and the measures involved in the two formulas
satisfy
$$\s_a = |E|^2\mu_\a.$$

We will call two $\Pi$-finite positive measures $\mu$ and $\ti\mu$ on $\R$  (AC) dual if there exists a Schur function $\varphi$ such that  $\mu=\s^\varphi_{-1}$ and $\ti\mu=\s^\varphi_1$, see Section \ref{periodic} and \cite{etudes} for detailed further discussion.

Connection between duality of measures, as defined above, to spectral problems is provided by the property that if $\mu$ is a spectral measure
of a CS with the Neumann initial condition at $0$, then   the spectral measure $\ti\mu$ of the same system with the Dirichlet initial condition is dual to $\mu$. Similarly, 
measures satisfying \eqref{eqClark} with other $\a$ correspond to other initial conditions for the same system. Note also that any two measures
$\s_\a$ and $\s_\b$, $\a\neq\b$, from the same AC family are dual up to constant
multiples, i. e., there exist $C_1,C_2>0$ such that $C_1\s_\a$ and $C_2\s_\b$ are dual.

\subsection{PW systems and PW measures}\label{PW}
Let $\HH$ be a Hamiltonian of a canonical system \eqref{eq001A}. We say that $\HH$ is of PW type ($\HH\in \PW$) if for any $t > 0$ there exists $s = s(t) > 0$ such that $s(t) \to \infty$ as $t\to\infty$ and the spaces $\BB_t(\HH)$ and ${\rm PW}_s$ are equal as sets (with possibly different norms), 
\begin{equation} \label{eq20}
    \BB_t(\HH) \doteq {\rm PW}_s. 
\end{equation}
We call the corresponding system a PW type system.

All regular Dirac systems (systems with locally summable potentials) are of PW type as follows from Lemma 3.1 in \cite{etudes} and the results of \cite{LS}; see also \cite{BR}.
The class of all PW type systems is significantly broader than the class considered in the classical Gelfand-Levitan theory; see \cite{etudes} for further discussion.

By definition, a positive measure $\mu$ on $\R$ is PW sampling  ($\mu\in\PW$) if it is sampling for all Paley-Wiener spaces PW$_t$:
$$\forall t \quad \exists C > 0, \qquad \forall f \in {\rm PW}_t, \qquad C^{-1} \|f\| \le \|f\|_{L^2(\mu)} \le C\|f\|.$$

Given $\mu$ and $\delta > 0$ we say that an interval $l \subset \R$ is $\delta$-{\it massive} with respect to $\mu$ if
$$\mu(l) \ge \delta \textrm{ and } |l| \ge \delta.$$

The $\delta$-{\it capacity} of an interval $I\subset\R$ with respect to $\mu$, denoted by $C_\delta(I)$, is the maximal number of disjoint $\delta$-massive intervals intersecting $I$.

The set of all PW sampling measures admits the following elementary description
\begin{theorem}[\cite{etudes}]\label{PWmu} $\mu$ is  PW sampling if and only if 
\begin{enumerate}[label=(\roman*)]
    \item For any $x \in \R$, $\mu(x,x+1) \le \const$;
    \item For any $t > 0$ there exist $L$ and $\d$ such that for all $I,\ |I| \ge L$, $$C_\delta(I) \ge t|I|.$$
\end{enumerate}
\end{theorem}

As the most basic example, any measure of the form $\rho(x) m(x)$, where $0 < c < \rho(x) < C < \infty$ is PW sampling. For further examples see \cite{etudes}.

We say that a measure $\mu$ has locally infinite support if $\text{supp\ } \mu \cap [-C,C]$ is an infinite set for some $C > 0$, or equivalently if $\text{supp\ }\mu$ has a finite accumulation point.
For a periodic measure one can easily deduce the following
\begin{corollary}[\cite{etudes}] \label{PWper}
A positive locally finite periodic measure is a Paley-Wiener measure if and only if it has locally infinite support.
\end{corollary}

PW systems and PW sampling measures are related via the following statement. For $\mu \in $PW we denote by ${\rm PW}_t(\mu)$ the Paley-Wiener spaces $PW_t $ endowed with the equivalent norm from $L^2(\mu)$.

\begin{proposition}[\cite{etudes}]
Suppose $\det(\HH)=1$ a.e., and let $\mu$ be the (unique) spectral measure of the corresponding CS \eqref{eq001A}. Then
	$$\mu \in {\rm(PW)} \qquad \Leftrightarrow \qquad \HH \in {\rm(PW)}.$$
Moreover, if either holds, then
	$$\forall t,\quad \BB_t(\HH) \doteq {\rm PW}_t(\mu).$$
\end{proposition}

\subsection{Det-normalization}
We have the following
\begin{theorem}[\cite{etudes}] \label{t0001} 
If $\HH$ is of PW type, then
\begin{equation}
    \det \HH \ne 0 \quad {\rm a.e.}, \qquad \int_0^\infty\sqrt{\det\HH(t)}dt = \infty. \label{eq005A}
\end{equation}
\end{theorem}

The change of time $t \mapsto s = s(t)$ in the above theorem allows us to transform any PW type system, or more generally any canonical system satisfying \eqref{eq005A} into a canonical system with $$\det \HH=1\quad{\rm a.e.}$$ Details on this change of variable can be found in \cite{etudes}.
We will call such systems {\it det-normalized}. A regular (locally summable) Hamiltonian will remain regular under det-normalization.
	
From here until the end of Section \ref{Hom}, all systems considered are assumed to be det-normalized; results in Sections \ref{PrelimUncertainty}-\ref{uncertainty} on the \Sch equation do not necessarily satisfy this assumption.
	
\begin{remark} $    $
\begin{itemize}
    \item All regular Dirac systems are of PW type as follows from the results of \cite{etudes, LS, BR}.
    \item  As mentioned in \cite{etudes}, characterization of PW Hamiltonians is a difficult problem even in the diagonal case. So far we know that
    $$\det \HH \ne 0 \qquad \not \Rightarrow \qquad  PW$$
    but 
    $$\det\HH \ne 0, \quad \HH \in W^{1,1}_{\rm loc} \qquad \Rightarrow \qquad  PW,$$
    as follows from the property that regular Dirac systems correspond to CS with $W^{1,1}_{\rm loc}$-Hamiltonians.
    \item An example from \cite{etudes} shows 
    $$\HH=\begin{pmatrix} h_{11} & h_{12} \\ h_{12} & h_{22} \end{pmatrix} \in \PW\ \not\Rightarrow \
    \begin{pmatrix} h_{22} & -h_{12} \\ -h_{12} & h_{11} \end{pmatrix} \in \PW.$$
    
    The second Hamiltonian corresponds to the switch between the Dirichlet and Neumann initial conditions.
\end{itemize}
\end{remark}

\section{Inverse spectral problem} \label{Etudes}

\subsection{Inverse problem: recovery of $h_{11}=h^{\mu}$}
For  a PW system with the spectral measure $\mu$  the leading term of the Hamiltonian $h_{11}=h^{\mu}$ can be recovered from the reproducing kernels $k_t = K_t(z, 0) \in {\rm PW}_t(\mu)$ defined in Section \ref{CSdB}:
	
\begin{theorem}[\cite{etudes}]\label{t3}
Let $\mu \in \PW$ be the spectral measure of a system \eqref{eq001A} with the Hamiltonian $\HH$. Then $t\mapsto k_t(0)$ is an absolutely continuous function and 
\begin{equation}
    h^\mu(t) := \pi \frac d{dt} k_t(0).\label{eq003A}
\end{equation}
\end{theorem}

Systems with even spectral measures have diagonal Hamiltonians \cite{etudes}. If, in addition, the system is det-normalized, then $h_{22}=1/h^{\mu}$ and the inverse spectral problem is solved with the recovery of $h^{\mu}$. For non-diagonal cases further analysis must be conducted, see below.

\subsection{Generalized Hilbert transform}\label{secHT}
For a $\Pi$-finite measure $\mu$ and $f \in L^2(|\mu|)$ we will use the notation
$$K(f\mu)(z) = \frac1{\pi} \int \frac{f(s)~d\mu(s)}{s - z},$$
and
$$\KK \mu(z) = \frac1{\pi} \int \left[ \frac1{s - z} - \frac{s}{1 + s^2} \right]~d\mu(s),$$
where $z \in \C \setminus \R$.
If $f \in L^2(\mu)$ is an entire function, then we define
$$H^\mu f = K(f\mu) - f \KK\mu.$$
It is clear that $H^\mu f$ extends to an entire function:
$$(H^\mu f)(z) = \frac1\pi \int \left[ \frac{f(s)-f(z)}{s-z} + \frac{s f(z)}{1+s^2} \right]~d\mu(s).$$

As was shown in \cite{etudes}, $H^\mu$ plays an important role in the recovery of the off-diagonal elements of the Hamiltonian.

\subsection{Inverse problem: recovery of $h_{12} = h_{21}$ ($g^\mu$)}
To recover the off-diagonal terms we use
\begin{theorem}[\cite{etudes}] \label{t007}
Let $\mu \in \PW$. Consider the reproducing kernels at 0, $k_t \in \PW_t(\mu)$.
Define $\tilde l_t = H^\mu k_t$.
Then $\mu$ is the spectral measure of the Hamiltonian
$$\HH = \begin{pmatrix} h^\mu & g^\mu \\ g^\mu & \frac{1+g^2_{\mu}}{h^\mu} \end{pmatrix},$$
where
$$g^\mu(t) := \pi\frac d{dt} \tilde l_t(0).$$
\end{theorem}

\subsection{Equations for the Fourier transform of $k_t$}\label{secFt} 
How to compute the functions $h^\mu$ and $g^\mu$ from Theorem \ref{t007}?
Denote by $k_t,\ \mathring k_t$ for the reprokernels $K_0,\ \mathring K_0$ in $\BB_t$ and PW$_t$ respectively. Sometimes it is helpful to work with the functions
$$\psi_t := \hat k_t,$$
so that
$$\frac1{\sqrt{2\pi}}\int_\R\psi_t = k_t(0).$$


If $f \in {\rm PW}_t$, then
$$f(0) = \frac1{\sqrt{2\pi}} \int_{-t}^t \hat f(\xi)~d\xi$$
while  also
$$f(0) =(f, \mathring k_t)_{PW_t}=\left(\hat f,\FF(\mathring k_t)\right)_{L^2(-t,t)}.$$

This leads us to the well known formula for the Fourier transform of the sinc function:
$$\FF \mathring k_t =\FF\left(\frac{\sin tz}{\pi z}\right)
=  \frac1{\sqrt{2\pi}}~1_{(-t,t)}.$$

The Fourier transform $\psi_t$ of the reprokernel $k_t$ can be found using the following statement.

\begin{theorem}[\cite{etudes}] \label{tTT}
$\psi=\psi_t$ satisfies \begin{equation}
    \psi \ast \hat\mu = 1 \quad  {\rm on} \quad (-t,t) \label{eq002A}
\end{equation}
and
$$\psi = 0 \quad  {\rm on}\quad \R \setminus [-t,t].$$
\end{theorem}
	
In general, all Fourier transforms and convolutions throughout the rest of the paper are understood in the sense of distributions.

\section{Periodic spectral measures and periodic approximations}\label{periodic}

\subsection{Measures from $M^+(\T)$}\label{secPer} 
Let $M^+(\T)$ be the set of all finite positive measures on the unit circle $\T \sim [0,2\pi]$ with {\it infinite support}. 
We identify measures on the circle with $2\pi$-periodic measures on $\R$, $M^+(\T) \subset M^+(\R)$, and define $M^+_{\rm even}(\T)\subset M^+(\R)$ as the subset of even periodic measures.

Notice that all measures in $M^+(\T)$ are $\Pi$-finite. Furthermore, as was mentioned before, it follows from Theorem~\ref{PWmu} (or Corollary~\ref{PWper}) that
\[
M^+(\T)\subset \PW.
\]
Let us list two more important properties of $M^+(\T)$
pertaining to our problems.

\begin{prop}[\cite{etudes}]
If $\mu\in M^+(\T)$, then all dual measures are in $M^+(\T)$.
\end{prop}

\begin{prop}[\cite{etudes}]
    If $\mu \in M^+_{\rm even}(\T)$, then there is a unique $\tilde \mu \in M^+_{\rm even}(\T)$ such that $\mu$ and $\tilde \mu$ are dual:
    $\ti\mu = \ti\mu_b,\ b = -\Re\KK\mu(0)$.
\end{prop}

Note that for the unique even dual measure $\ti\mu$ we have
\begin{equation}
    \KK\ti\mu = -\frac 1{\KK\mu}, \ \PP\ti\mu = -\Re \frac 1{\KK\mu}.   \label{eqEven}
\end{equation}

\subsection{Moments}\label{secMom} 
If $\mu \in M^+(\T)$, we can consider the Fourier series
  $$\mu \sim \sum_{-\infty}^\infty \gamma_k e^{ikx},$$
where $$\gamma_k = \frac1{2\pi} \int_0^{2\pi} e^{-ikx}~d\mu(x) = \a_k + i\b_k$$ are the "trigonometric moments" of $\mu$. 
Clearly, $$\gamma_{-k} = \overline{\gamma_k},$$
so $$\mu\sim \gamma_0 + 2\sum_1^{\infty}(\alpha_k\cos kx - \beta_k\sin kx),$$
and $\mu$ is even iff all moments are real, i.e., all $\beta_k = 0$.

The moments can be arranged into the Toeplitz matrix
$$\Gamma(\mu) = \begin{pmatrix} \gamma_0 & \gamma_1 & \gamma_2 & \dots \\
  \gamma_{-1} & \gamma_0 & \gamma_1 & \dots\\
  \gamma_{-2} & \gamma_{-1} & \gamma_0 & \dots\\
  \dots & \dots & \dots & \dots
  \end{pmatrix}.$$
We will denote by $\Gamma_n(\mu)$ the  $n \times n$  matrix in the upper left corner of $\Gamma$.
Note that our matrix index starts from 0 instead of 1.
The well-known Caratheodory theorem says that a sequence $\{\gamma_k\}$ is the sequence of moments of some $\mu \in M^+(\T)$ iff $$\forall n, \qquad \det \Gamma_n > 0.$$

Finally let us recall the relation between the moments of a periodic measure and its Fourier transform: if $\mu \sim \sum \gamma_k e^{ikx}$, then
$$\hat\mu = \sqrt{2\pi} \sum_{-\infty}^\infty \gamma_k\delta_k.$$
  
\subsection{Computation of $h^\mu(t)$}\label{secHmu}
Recall that if $\mu\in\PW$, then
    $$h^\mu(t):=\pi\frac d{dt}\left[L_{\mu,t}^{-1} \stackrel \circ{k_t}(0)\right]$$

If $A = (a_{jk})$ is a matrix we use the notation $\SS[A]$ for the sum of the elements of $A$:
    $$\SS[A]=\sum_{ j,k} a_{jk}.$$

\begin{theorem}[\cite{etudes}]\label{tPer} 
If $\mu \in M^+(\T)$, then the function $h^\mu(t)$ is locally constant on $\R_+\setminus \frac12\mathbb N$, i.e., 
$$h^\mu(t)=h_0,\; h_1, \; h_2,\; \dots \quad {\rm on} \quad \left[0, \frac12 \right), \; \left[ \frac12, 1\right ),\; \left[ 1, \frac32 \right), \; \dots,$$
and
$$h_n = \SS\left[\Gamma_{n}(\mu)^{-1}\right] - \SS\left[\Gamma_{n - 1}(\mu)^{-1}\right], \quad n \geq 1, \quad h_0 = \SS[\Gamma_0^{-1}] = \gamma_0^{-1}.$$
\end{theorem}

Via a change of variable this result holds for $2T-$periodic measures with $h^\mu(t)$ being locally constant on $\R_+ \setminus \frac{\pi}{2T}\mathbb N$.

This explicit inverse algorithm can be reversed to solve the direct spectral problem. We have a finite-dimensional algorithm as well.

\begin{theorem}[\cite{Direct}]\label{directalg}
Let $\mathcal{H}$ be a diagonal Hamiltonian that is locally constant on $\R_+\setminus \frac12\mathbb N$, i.e., 
$$h^\mu(t)=h_0,\; h_1, \; h_2,\; \dots \quad {\rm on} \quad \left[0, \frac12 \right), \; \left[ \frac12, 1\right ),\; \left[ 1, \frac32 \right), \; \dots,$$
where all $h_n$ are positive real numbers.
Then the spectral measure $\mu$ is a $2\pi-$periodic measure with moments
\begin{equation*}
\begin{aligned}
    \gamma_0 &= \frac{1}{h_0}, \quad \gamma_1 = \frac{1}{h_0} \cdot\frac{1 - h_1}{1 + h_1}, \\
    \gamma_n &= \frac{1}{h_{11}^0} \cdot \frac{\left( 1 + D_n/\Delta_n \right) \left( 1 - \mathds{1}\Gamma_{n - 1}^{-1} \mathbf{u}_n \right)^2 - (\Delta_n - D_n) h_{n + 1}}{h_{n + 1} \Delta_n + \left( 1 - \mathds{1}\Gamma_{n - 1}^{-1} \mathbf{u}_n \right)^2/\Delta_n},
\end{aligned}
\end{equation*}
where $\mathbf{u}_n = (c_1, c_2, \ldots, c_n)^T$, $\mathbf{v}_n = (c_n, c_{n - 1}, \ldots, c_1)^T$, $\mathds{1} = (1, 1, \ldots, 1)$, $\Delta_n = 1 - \mathbf{u}_n^T \Gamma_{n - 1}^{-1} \mathbf{u}_n$ and $D_n = \mathbf{u}_n^T \Gamma_{n - 1}^{-1} \mathbf{v}_n$.
\end{theorem}

\begin{remark}
When $\gamma_0 = 1$, the algorithm simplifies significantly, since the Verblunsky coefficients $\alpha_n$ of the spectral measure $\mu$ are then well-defined. The Verblunsky coefficients of $\mu$ can be obtained using
\[
    \alpha_n = \frac{1 - h_{n+1}/h_n}{1 + h_{n+1}/h_n}.
\]

Also note that there is no growth rate or decay rate restriction on the sequence $\{h_n\}$. In Section \ref{examples}, we provide the explicit solution for $h_n = a^n$, $a > 0$.
\end{remark}

\subsection{Computation of $g^\mu(t)$}\label{secGmu1} 
Let $\mu\in M^+(\T)$ and let
$$\Delta(\mu) = \begin{pmatrix} 0 & \gamma_1 & \gamma_2 & \dots\\
  -\gamma_{-1} & 0 & \gamma_1 & \dots\\
  -\gamma_{-2} & -\gamma_{-1} & 0 & \dots\\
  \dots & \dots & \dots & \dots
  \end{pmatrix}$$
where, as before, $\gamma_k$ are trigonometric moments of $\mu$.
We denote by $\Delta_n$ the $n\times n$ matrix in the upper left corner of $\Delta$.

\begin{theorem}[\cite{etudes}]\label{tPer2}
If $\mu\in M^+(\T)$, then the function $g^\mu(t)$ is locally constant on $\R_+\setminus \frac12\mathbb N$, i.e. $$g^\mu(t) = g_0,\; g_1, \; g_2,\; \dots \quad {\rm on}\quad \left[0,\frac12\right),\; \left[\frac12,1\right),\; \left[1,\frac32\right),\;\dots$$ and
$$g_n = \SS\left[\Delta_{n}\Gamma_{n}^{-1}\right]-\SS\left[\Delta_{n - 1}\Gamma_{n - 1}^{-1}\right], \quad n \geq 1, \quad g_0 = 0.$$ 
\end{theorem}

\subsection{Canonical systems and orthogonal polynomials on the unit circle}
\label{secOPUC}
In this section we investigate the connection between the inverse spectral problems for CS with periodic spectral measures and orthogonal polynomials on the unit circle.
We show that the formula for $h_n$ from Theorem \ref{tPer} and the formula for $g_n$ from Theorem \ref{tPer2} can be interpreted as point evaluations of such polynomials.

For $\mu\in M_+(\T)$ we denote by $\ONP(\mu)$ the family of polynomials of $z$ on the unit circle orthonormal in $L^2(\mu)$.
By $\phi_n \in \ONP(\mu)$ we denote the polynomial of degree n.

\begin{theorem}[\cite{etudes}]\label{tONP}
Let $\mu \in M^+(\T)$ and let $\{\phi_n\}$ be \ONP($\mu$). We have
    $$ h^\mu_n = |\phi_n(1)|^2.$$
\end{theorem}

Theorem \ref{tONP} allows us to describe the change of $h_n$ corresponding to a shift of a periodic spectral measure.
For $\eta\in\mathbb T$ we define $\mu^{(\eta)}(\zeta)=\mu(\eta\zeta)$.

\begin{corollary}[\cite{etudes}]
    $$ h_n^{(\eta)} = |\phi_n(\eta)|^2.$$
\end{corollary}

For the off-diagonal term $g^\mu_n$, we have a similar result. We recall from Section \ref{secClark}, $\tilde \mu $ is AC dual to $\mu$ if there exists a function $\phi \in H^\infty(\C_+)$, $\|\phi\|_\infty \le 1,$ such that
$$\mu = \sigma_{-1}^\phi, \qquad \tilde\mu = \sigma_{1}^\phi,$$
i.e.,
\begin{equation}
\PP\mu = \Re\frac{1-\phi}{1+\phi}, \qquad \PP\tilde\mu = \Re\frac{1+\phi}{1-\phi}.\label{eq010}
\end{equation}

An alternative way to define the AC dual is through the Verblunsky coefficients of $\mu$: let $\{\alpha_n\}$ be the Verblunsky coefficients of the measure $\mu \in M^+(\T)$, the AC dual $\tilde \mu$ is the measure in $M^+(\T)$ satisfying \begin{equation*}
    \alpha_n(\tilde \mu) = - \alpha_n.
\end{equation*}
Proof of equivalence of these two definitions can be found for instance in \cite{Simon}.

\begin{theorem}[\cite{h12}]\label{tONP2}
Let $\mu \in M^+(\T)$, and denote by $\tilde \mu$ the Clark dual of $\mu$. Let $\{\phi_n\}$ be \ONP($\mu$) and $\{\Tilde{\phi}_n\}$ be \ONP($\tilde \mu$). We have
    $$g^\mu_n = -\frac{\Im \left({\Tilde{\phi}_n(1)/\phi_n(1)}\right) }{\Re \left({\Tilde{\phi}_n(1)/\phi_n(1)}\right) }.$$
\end{theorem}

\subsection{Periodic approximations}
Until recently, examples of solutions to ISP for canonical systems were rare in the literature. 
Even though the number of available examples has increased with the development of new methods, 
most of the new examples pertain to periodic spectral measures. 
Non-periodic examples are still difficult to calculate explicitly, even for the simple spectral measures. 
The results of this section allow one to find such examples numerically after using the theorems of the previous section to find the Hamiltonians of the periodizations of the spectral measure.

Let $\mu$ be a measure on $\R$. We denote by $\mu_T$ its $2T-$periodization, defined as $\mu_T(S) = \mu(S)$ for $S \subset [-T, T)$, and  extended periodically to the rest of $\R$.

By Corollary \ref{PWper}, a locally finite periodic measure is a PW measure if and only if it has infinite support on its period. If a locally finite $\mu$ has infinite support on an interval $[-C, C]$ for some $C > 0$, then $\mu_T$ is a PW measure for $T > C$. Therefore, we can recover a PW Hamiltonian corresponding to $\mu_T$ using theorem \ref{tPer} or \ref{tONP}. We denote by $h^{\mu}_T$ the upper-left entries of the Hamiltonians corresponding to $\mu_T$.

In various classes of spectral measures $\mu$, including some not in the PW class, we can show that $h^{\mu}_T$ converges to $h^{\mu}$ as $T \to \infty$.

Recall that $K_0^t$ are the reproducing kernels in the de Branges' spaces $B_t(\mu)$. Denote by $K_0^{t, T}$ the reproducing kernels in the de Branges' spaces $B_t(\mu_T)$. We begin with the following lemma on the convergence of $h^{\mu}_T$ to $h^{\mu}$.

\begin{lemma}[\cite{PZ}]\label{SuffCondition}
If $ \lim_{T \to \infty} K_0^{t, T}(0) = K_0^t(0)$ for almost every $t$, then for any
interval $(a,b)\subset \R_+$, 
$$\int_a^b h^{\mu}_T(t)dt  \to \int_a^b h^{\mu}(t)dt $$  as $T \to \infty$.
\end{lemma}

Since the Hamiltonian matrix $\mathcal{H}(t) \geq 0$ almost everywhere, $h^{\mu} \geq 0$ almost everywhere.
For non-negative functions convergence of integrals over intervals, like in the last lemma, 
is equivalent to the convergence on continuous compactly supported functions

\begin{corollary}[\cite{PZ}]
If $ \lim_{T \to \infty} K_0^{t, T}(0) = K_0^t(0)$ for almost every $t$, then
$$\int_0^\infty h^{\mu}_T(t)\phi(t)dt{\to} \int_0^\infty h^{\mu}(t)\phi(t)dt$$
as $T\to\infty$, for every continuous compactly supported function $\phi$ on $\R_+$.
\end{corollary}

Let  $f_T, T > 0$ and $f$ be $L^1_{loc}(\R_+)$-functions.
We write $f_T \overset{*}{\to} f$ as $T \to \infty$ if $f_T, T>0$ and $f$ satisfy the conclusions of the last lemma and corollary in place of $h^{\mu}_T$ and $h^{\mu}$. 
Using this notation, $ \lim_{T \to \infty} K_0^{t, T}(0) = K_0^t(0)$ for almost every $t$ implies $h^{\mu}_T\overset{*}{\to} h^{\mu}$ as $T\to \infty$. 

First, we study convergence of periodizations in the class of canonical systems whose spectral measures are from the PW class. 
As was mentioned before, it extends the class of regular systems considered in classical literature and was recently studied in \cite{BR, Bessonov, etudes}.


For PW systems we prove the following convergence result, which implies Hamiltonians of PW systems can be recovered by Hamiltonians of their periodizations. 

\begin{theorem}[\cite{PZ}]\label{PWConv}
Let $\mu$ be an even PW measure with locally infinite support. Then
\begin{equation*}
    h^{\mu}_T\overset{\ast}{\to} h^{\mu}\text{ as }T\to\infty.
\end{equation*}
\end{theorem}

Motivated by Example \ref{pointmassexample} in \cite{etudes}, we extend this result to show that the periodization approach can, in some cases, be applied beyond the PW class.

It follows from theorem \ref{PWmu} that any measure decaying to $0$ at infinity is not a PW measure. Nonetheless, we have the following statement, which shows that Hamiltonians of periodizations can be used to solve inverse problems with decaying spectral measures.

\begin{theorem}[\cite{PZ}]\label{decay}
Let $\mu$ be a positive even Poisson-finite measure on $\R$. Suppose that there exist $C,c > 0$, such that $\mu$ has infinite support on $[-C, C]$, and $\mu(x, x + c)$ is decreasing on $\mathbb{R}_+$. Then, for $T_n=nc$, $$h_{T_n}^{\mu} \overset{\ast}{\to} h^{\mu}\text{ as }n\to\infty,\ n\in\mathbb N.$$
\end{theorem}

The following particular case can be proved in a similar way or deduced from the last statement.

\begin{corollary}[\cite{PZ}]\label{CorDecay}
Let $\mu$ be a positive even absolutely continuous locally-finite measure, whose density decreases on $\R_+$. Then \begin{equation*}
    h^\mu_T \overset{\ast}{\to} h^{\mu} \text{ as } T \to \infty.
\end{equation*}
\end{corollary}

Our next theorem allows one to further increase the class of measures for which the periodization approach works.

\begin{theorem}[\cite{PZ}]\label{polygrowth}
Let $\mu$ be an even Poisson-finite measure such that $d\mu = h(x) d\nu$, where $\nu$ satisfies the conditions of theorem \ref{decay}, and $h(0) \neq 0$ is an even positive function increasing on $\mathbb{R}_+$, and satisfies $h(x) = O(|x|^q)$ as $|x|\to\infty$ for some $q>0$. Then, for $T_n=nc$,
$$h_{T_n}^{\mu} \overset{\ast}{\to} h^{\mu}\text{ as }n\to\infty,\ n\in\mathbb N.$$
\end{theorem}

The following particular case can be proved in a similar way or deduced from the last statement.

\begin{corollary}[\cite{PZ}]\label{CorPolygrowth}
Let $\mu$ be a positive even absolutely continuous locally-finite measure, whose density is the product of $h(x)$, an even positive function increasing on $\mathbb{R}_+$, satisfying $h(x) = O(|x|^q)$ as $|x| \to \infty$ for some $q > 0$, and $f(x)$, an even non-negative function decreasing on $\mathbb{R}_+$. Then \begin{equation*}
    h^{\mu}_T\overset{\ast}{\to} h^{\mu}\text{ as }T\to\infty.
\end{equation*}
\end{corollary}

\subsection{Step function approximations}
Theorem \ref{tPer} shows that if a diagonal $\HH$ is a step function of uniform step size, then the corresponding spectral measure $\mu$ is an even periodic PW measure. Moreover, we can solve this direct spectral problem with Theorem \ref{directalg}. 

If $\HH = \begin{pmatrix} h_{11} & 0 \\ 0 & h_{22}\end{pmatrix}$ is not a step function, we can 'periodize' the Hamiltonian by defining a step-function Hamiltonian $\HH^T$ as follows: \begin{equation}\label{periodization}
\begin{aligned}
    h_{11}^{T, n} &= h_{11}^T \big|_{\left[ nT, (n + 1)T \right)]} = \frac{1}{T} \int_{nT}^{(n + 1)T} h_{11}(s) ds,\\
    h_{22}^{T, n} &= h_{22}^T \big|_{\left[ nT, (n + 1)T \right)]} = \frac{1}{h_{11}^{T, n}}.
\end{aligned}
\end{equation}
The spectral measure corresponding to the canonical system with Hamiltonian $\HH^T$ is denoted by $\mu_T$.

Given that all $\mu_T$'s are PW measures, our focus is on applying this approach to PW systems, since all their de Branges spaces $B(E^T_t)$ are the same as $PW_t$ as sets. Unfortunately, characterizing PW systems solely through their Hamiltonians is challenging as discussed earlier in Section \ref{prelim}. Describing PW Hamiltonians that are not step functions, even in the diagonal case, is inherently difficult. It is known that not all Hamiltonians with non-vanishing determinants are PW Hamiltonians. However, PW-Hamiltonians include those arising from real Dirac systems, as shown in several results in \cite{etudes, BR}.

Real Dirac systems on the half-line $\mathbb{R}_+$ are of the form \begin{equation}\label{RD}
    \Omega \Dot{X}(t) = z X(t) - Q(t) X(t), \ \ t \in (0, \infty).
\end{equation}
Here, $z\in \mathbb{C}$ is a spectral parameter,
$\Omega = \begin{pmatrix} 0 & 1 \\ -1 & 0 \end{pmatrix}$ is the same symplectic matrix as in canonical systems, and $Q(t) = \begin{pmatrix} 0 & f(t) \\ f(t) & 0 \end{pmatrix}$ is the potential matrix, where $f(t)$ is a real-valued locally integrable function. These systems, when rewritten as canonical systems, feature det-normalized diagonal Hamiltonians, with
$$h_{11}(t) = \exp \left( \int_0^t f(s) ds \right),$$
where $f$ is the locally integrable function in $Q(t)$. Details on rewriting such a system as a canonical system can be found, for instance, in \cite{MPS, Rom, Zhang}.

When using the periodization approach on Hamiltonians from real Dirac systems, $\HH^T$'s will still be in the form of \begin{equation*}
    \begin{pmatrix}
        \exp \left( \int_0^t f^T(s) ds \right) & 0 \\ 0 & \exp \left( - \int_0^t f^T(s) ds \right)
    \end{pmatrix},
\end{equation*}
but with $f^T$ being a discrete measure supported on an arithmetic progression.

\begin{theorem}[\cite{Direct}]\label{directconv}
Let $\HH$ be the Hamiltonian of a real Dirac system, and let $H^T$'s be the step-function approximations defined in \eqref{periodization}. Denote by $\mu$ the spectral measure of the real Dirac system, and $\mu_T$ the spectral measure of the canonical system with Hamiltonian $\HH^T$. Then, for $\phi \in PW_a$, $a > 0$, we have $\|\phi\|_{L^2(\mu_T)} \to \|\phi\|_{L^2(\mu)}$ as $T \to 0$.
\end{theorem}

As noted at the beginning of this section, the results here, together with those from Section \ref{Etudes}, can be used to construct explicit examples of spectral problems. We present a variety of such examples at the end of the paper in Section \ref{examples}.

\section{Homogeneous measures and spaces} \label{Hom}

\subsection{Homogeneous spectral measure}\label{Hommea}
For a measure $\mu(x)$ on $\R$ and $t > 0$ we denote by $\mu_t(x)=\frac 1t\mu(tx)$ the measure such that for any Borel $B\subset \R$,
$$\mu_t(B)=\frac 1t\mu(tB).$$

A measure $\mu$ is homogeneous if 
$$\forall t > 0,\qquad \mu_t(x) = \mu(x).$$

It is easy to show that a homogeneous measure must be absolutely continuous, $d\mu(x)=\rho(x)dx$, where $\rho(x)$ is constant on $\R_+$ and on $\R_-$.
We will assume that the constants are strictly positive, so that $\mu \in \text{PW}$.

Note that all examples of such measures measures are given 
by \begin{equation}
\mu = c_1m + c_2\sigma, \label{eq005}
\end{equation}
where $\sigma(x) = \sign(x)$, and $c_1>|c_2|$, $c_1,c_2\in \mathbb R$.

Also, from the well-known rescaling properties of the Hamiltonians and spectral measure, it follows that $h^{\mu}(t)$ and $g^\mu(t)$ are both constant functions. The constant $h^\mu(t)$ can be obtained via direct calculations \cite{etudes2}, and a sketch of the calculation is presented in Section \ref{ISPHom}.

To find $g^\mu(t)$, we need the following results. As before we denote by $k_t(z)$ the reprokernels of the dB spaces $\PW_t(\mu)$ at zero.
\begin{theorem}[\cite{etudes2}]\label{homid}
If $\mu\in \PW$ is homogeneous, then we have the identity
$$k_t(z)=tk_1(tz).$$
\end{theorem}

\begin{example}
If $d\mu=dx$, then 
$$k_t(z)\equiv k_t^0(z)=\frac1\pi~\frac {\sin tz}z,$$
and 
$$k_t^0(z)=tk_1^0(tz).$$ 
\end{example}
 
Following \cite{dB},  we say that a dB space $\BB=\BB(E)$ is homogeneous   if for all $t\in(0,1)$
$$F\in \BB\quad\Rightarrow\quad  t ^{1/2} F(tz)\in\BB$$
and both functions have the same norm in $\BB$.

\begin{theorem}[\cite{etudes2}]\label{thm2}  
$\mu \in \PW$ is homogenous iff all its dB spaces are homogenous. 
\end{theorem}

\begin{corollary}[\cite{etudes2}]\label{cor} Let $d\mu=\rho(x)dx$.
The following  three conditions are equivalent:
\begin{enumerate}[label=(\roman*)]
    \item $\rho(x) = \rho(tx)$ for all $t \in (0,1)$;
    \item $\BB_t$ is homogeneous  for all $t > 0$;
    \item $k_t(z) = tk_1(tz)$  for all $t\in (0,1)$.
\end{enumerate}
\end{corollary}

\subsection{Quasi-homogeneuous spectral measures} 
By definition, a measure $\mu(x)$ is quasi-homogeneous of order $\nu$ if
for all $t > 0$,
$$t^{1+2\nu} \mu(x) = \mu_t(x).$$
Once again, it is not difficult to prove that such measures are absolutely continuous. 
Their densities $\rho$ must satisfy
$$\forall t > 0,\quad t^{1+2\nu} \rho(x/t) = \rho(x),$$
or equivalently 
$$\forall t > 0,\quad t^{1+2\nu} \rho(y) = \rho(ty).$$

One can show that quasi-homogeneous measures form a two parameter family:
$$\rho(x) = \begin{cases}
x^{1+2\nu} \rho(1),\quad x > 0\\
|x|^{1+2\nu} \rho(-1),\quad x < 0
\end{cases},$$
with arbitrary positive constants $\rho(\pm 1)$.
As usual, a special case among spectral measures of canonical systems is occupied by even measures, which correspond to the case $\rho(1)=\rho(-1)$. 

The measures are not Paley-Wiener unless $\nu = -\frac{1}{2}$.
Quasi-homogeneous measures are locally finite (on $\R$) iff $\nu>-1$, and
Poisson-finite iff $\nu<0$.
For this reason, we will only consider the case
\[
-1<\nu<0.
\]
The value $\nu=-1/2$ was the homogeneous case discussed in the previous subsection.



The results of this section can be applied to obtain explicit solutions of ISP for a homogeneous system in the PW case and to relate spectral problems to Bessel functions. These applications are discussed in detail in \cite{etudes2}. A compact presentation of these ideas is given at the end of the paper in Section \ref{examples}, with the ISP for a homogeneous system in Section \ref{ISPHom} and the connection with Bessel functions in Section \ref{bessel}.

We now shift gears slightly to discuss mixed spectral problems, originating in the work of Hochstadt and Lieberman \cite{HL} and further developed by del Rio, Gesztesy, and Simon \cite{DGS}. Before doing so, we briefly introduce some additional background needed for this topic.

\section{Holomorphic functions in spectral problems} \label{PrelimUncertainty}
\subsection{The \Sch operator and its spectra}\label{spectra}

We consider a special case of CS, the \Sch equation \begin{equation}
Lu = -u'' + qu = z^2 u \label{e00}
\end{equation}
on the interval $[0,\pi]$. To avoid unnecessary technicalities we will make the following assumptions.

Even though some of our results can be extended to the case $q \in L^p$, we mostly restrict ourselves to square summable potentials,  $q \in L^2$, which makes some of the statements considerably shorter.
For the same reason we will consider only Dirichlet ($u = 0$) or Neumann ($u' = 0$) boundary conditions at the endpoints.

The operator $$L : u \mapsto -u'' + qu $$ is a self-adjoint operator in $L^2([0,\pi])$ whose domain is the set of functions $u$ with absolutely continuous derivatives which satisfy the
boundary conditions and the condition $-u'' + qu \in L^2([0,\pi])$.
In the case $q \in L^2([0,\pi])$ the domain consists of functions from the Sobolev space $W^{2,2}$ satisfying the boundary conditions.
The spectrum of the operator with $q \in L^1$ and any self-adjoint boundary conditions at the endpoints is bounded from below.
We will assume the spectrum $\Sigma_{ND}$ of $L$ with Dirichlet-Neumann boundary conditions to be positive.
This assumption is not overly restrictive as any operator can be made positive by adding a positive constant to $q$: the spectrum then shifts to the right by the same constant.

We will denote the set of operators satisfying the above conditions by $\SS^2$. Occasionally we will use the notation $\SS^1$ for the operators with summable potentials.

As it is often done in this area, to facilitate the application of standard tools of Fourier analysis we will apply the square root transform to the spectra and analytic functions associated with our operators. 
In particular, we will denote by $\sigma_{DD}$ and $\sigma_{ND}$ the spectra of the operator $L$ after the square root transform:
$$\sigma_{DD} = \{ \lan | \lan^2 \in \Sigma_{DD} \} \cup \{ 0 \}, \quad \sigma_{ND} = \{ \lan | \lan^2 \in \Sigma_{ND} \}.$$
Note that under our positivity assumptions  both spectra are real.

We say that a sequence of complex numbers is discrete if it does not have finite accumulation points.
Both spectra $\sigma_{DD}$ and $\sigma_{ND}$ are real discrete sequences. Moreover, it is well known that
\begin{equation}
\sigma_{DD} = \{ \lan \}_{n \in \Z},\ \l_0 = 0, \ \l_{\pm n} = \pm \pi \left( n + \frac Cn + \frac {a_{n}}n \right),\ n \in \N, \label{e000a}
\end{equation}
and
\begin{equation}
\sigma_{ND} = \{ \eta_n \}_{n \in \pm \N},\ \eta_{\pm n} = \pm \pi \left( n - \frac 12 + \frac Cn + \frac {b_n}n \right),\ n \in \N ,\label{e000}
\end{equation}
where $C$ is a real constant and $\{a_n\},\{b_n\} \in l^2$.
Conversely, two interlacing sequences are equal to the spectra $\sigma_{DD},\ \sigma_{ND}$
of a \Sch operator on $[0,\pi]$ with an $L^2$-potential if and only if the sequences satisfy the above asymptotics with some $a_n,b_n\in l^2$.
These asymptotics follow from more general formulas by Marchenko \cite{M2}, see also \cite{Tru} or \cite{BBP}.

\subsection{Analytic integrals of spectral measures}\label{functions}

We denote by $L^1_\Pi$ the space of Poisson-summable functions on $\R$.
If $u \in L^1_\Pi$  we define its Herglotz integral as
$$ H u(z) = \frac1{\pi} \int \left[ \frac 1{t - z} - \frac t{1 + t^2} \right] u(t)dt.$$
If $\mu$ is a Poisson finite measure on $\R$, then
$$ H \mu(z) = \frac1{\pi} \int \left[ \frac 1{t - z} - \frac t{1 + t^2} \right] d\mu(t).$$

Denote by $u(t, z) $ the solution of \eqref{e00} with boundary conditions $u(0, z)=0,\ u'(0, z)=1$, where the derivative is taken on $t$. The Weyl function $m_+$ is defined as
 $$ m_+(z) = -\frac{u'(\pi, z)}{z u(\pi, z)}.$$
It is well-known that $m_+$ is a Herglotz integral of a positive measure supported on $\sigma_{DD}=\{\lan\}$: 
$$m_+(z) = H\mu_+(z),$$ 
$$\mu_+ = \alpha_0 \delta_0 + \sum_{n \in \N} \alpha_n(\delta_{\lambda_{n}} + \delta_{\l_{-n}}),$$
\begin{equation}
    \alpha_{n} = 1 + \frac{c_n}{n + 1}, \ c_n \in l^2.\label{asymp1}
\end{equation}

Similarly, if $v(t, z)$ is the solution of \eqref{e00} with boundary conditions $v(\pi, z) = 0,\ v'(\pi, z) = 1$, one defines $$m_- = \frac{v'(0, z)}{z v(0, z)}.$$
Then $\mu_-$ is defined via
$$m_-(z) = H \mu_-(z),$$
$$\mu_- = \beta_0 \delta_0 + \sum_{n \in \N} \beta_n(\delta_{\lambda_{ n}} + \delta_{\l_{-n}}),$$
\begin{equation}
    \beta_{ n} = 1 + \frac {d_n}{n + 1},\ d_n \in l^2.\label{asymp2}
\end{equation}

The asymptotics of the pointmasses $\alpha_n, \beta_n$ can be deduced from the spectral asymptotics of $\sigma_{DD}$ and $\sigma_{ND}$ discussed in the last section. 
Moreover, together the asymptotics for $\lan \in \sigma_{DD}$ ($\eta_n\in \sigma_{ND}$) and for $\alpha_n$ ($\beta_n$) give an if and only if condition for a positive measure to be a spectral measure $\mu_+$ ($\mu_-$) of a \Sch operator from $\SS^2$.
The famous theorem by Marchenko \cite{M1} says that a regular ($q\in L^1$) \Sch operator can be uniquely recovered from its spectral measure $\mu_+$ (or $\mu_-$).
Another classical result is the following two-spectra theorem by Borg.
Throughout the paper, we will say that a \Sch operator is uniquely determined by its (mixed) spectral data, if any other operator from the same class with the same data must have the same potential $q$ a.e. on $[0, \pi]$.

\begin{theorem}[\cite{Borg}]\label{Borg} 
A \Sch operator $L\in \SS^1$ is uniquely determined by its spectra $\sigma_{DD}$ and $\sigma_{ND}$.
No proper subset of $\sigma_{DD}\cup\sigma_{ND}$ has the same property.
\end{theorem}

As was mentioned in the introduction, we will present an 'uncerainty' version of this theorem in Section \ref{uncertainty}.

For a subinterval $(a, b)$ of $(0, \pi)$ recall from Section \ref{CSdB} that $M = M_{(a, b)}$ is the transfer matrix, 
Note that by the Wronskian identity, $\det M \equiv 1$ and $M(w_z(a),w'_z(a))^T = (w_z(b),w'_z(b))^T$ for any solution $w_z$ of \eqref{e00}.

\subsection{The even operator}\label{even}

We will call an operator $L$ even (with respect to the middle of the interval), if its potential satisfies $q(x)=q(\pi-x)$ for all $x\in(0,\pi/2)$. Such an operator is uniquely
defined by one of its spectra, say $\sigma_{DD}$, see for instance \cite{Tru}. If $\L=\{\lan\}$ is a sequence satisfying \eqref{e000a}, we denote by $\Gamma(\L)=\{\gamma_n\}$ the sequence of pointmasses
of the spectral measure $\mu_+=\sum \gamma_n\delta_{\lan}$ corresponding to the unique even operator with  $\sigma_{DD}=\{\lan\}$. Notice that in the even case $\mu_+=\mu_-$.

\begin{lemma}[\cite{MP3}]\label{leven} For any \Sch operator $L\in \SS^1$ the pointmasses of $\mu_\pm$ satisfy
$$\alpha_n\beta_n=\gamma_n^2.$$

\end{lemma}

\begin{remark}\normalfont The constants $\gamma_n$, the pointmasses of the spectral measure of the even operator, appear under different names in other problems
of analysis. As an example, let us point out the following connection with the classical problem of completeness of polynomials in weighted spaces.

We assume that all discrete real sequences are enumerated in natural increasing order.

We say that a sequence $\L=\{\lan\}$ has (two-sided) upper density $d$ if
$$\limsup_{A\to\infty}\frac{\#[\L\cap (-A,A)]}{2A}=d.$$
If $d=0$ we say that the sequence has zero density.
  A discrete sequence $\L=\{\lan\}$ is called \textit{balanced} if the limit
\begin{equation}
\lim_{N\to\infty}\sum_{|n|<N}\frac {\l_n}{1+\lan^2}
\label{balance}
\end{equation}
exists.

Let $\L=\{\lan\}$ be a balanced sequence of finite upper density. For each $n, \ \lan\in\L,$ put
$$p_n= \frac 12\left[\log(1+\l_n^2)+\sum_{n\neq k, \ \l_k\in\L}\log\frac{1+\l_k^2}{(\l_k-\lan)^2}\right],$$
where the sum is understood in the sense of principle value, i.e. as
$$\lim_{N\to\infty}\sum_{0<|n-k|<N}\log\frac{1+\l_k^2}{(\l_k-\lan)^2}.$$
We will call the sequence of such numbers $P(\L)=\{p_n\}$ the \textit{characteristic sequence} of $\L$.

Here is a sample of a statement on completeness of polynomials in terms of characteristic sequences.

\begin{theorem}\cite{Polya, CBMS}
Let $\mu$ be a finite  positive discrete measure supported on   $\R$ such that $L^1(\mu)$ contains polynomials.

\no Polynomials are not dense in $L^1(\mu)$  if and only if there exists a balanced zero density subsequence $\L=\{\lan\}\subset \supp \mu$ such that its characteristic sequence $P(\L)=\{p_n\}$ satisfies
$$\exp{p_n}=O(\mu(\{\lan\}))$$
as $|n|\to \infty$.
\end{theorem}

Similar statements can be formulated for families of exponential functions in place of polynomials.
In such statements zero density sequences are replaced with sequences of positive density, which
makes them closer related to the regular spectral problems considered in this note. For such results, along with the case of $L^p,\ p\neq 1,$ or Bernstein's spaces, see for instance \cite{Polya, CBMS}.
The case of a general measure can be reduced to the discrete case via some of the standard tools of completeness problems.
For further references and historic remarks see also
\cite{Lub, Meg}.

As the reader may have already guessed, the characteristic sequence $p_n$ is nothing else but the sequence of pointmasses of the even operator
$\gamma_n$ in the case when $\L$ is a spectral sequence, i.e., $P(\L)=\Gamma(\L)$. Since any discrete sequence is a spectral sequence for a suitably chosen Krein's canonical system,
one could formulate the last statement using pointmasses of the even operator instead of characteristic sequences.

An alternative way to define the constants $\gamma_n$ is
$$\gamma_n=1/|F'(\lan)|,$$
where
$$F(z)=\prod \left(1-\frac z{\lan} \right).$$

\end{remark}

\subsection{Hermite-Biehler functions for \Sch operators}\label{HBSthm}

Not every Hermite-Biehler function can be obtained from a \Sch equation in the way described in the last section. 
Characterization of such functions is important for applications, as we will illustrate in Section \ref{sHorvath}.
Such a characterization was recently obtained in \cite{BBP}. Here we
present a version of the same result which will be more convenient for our purposes.

Consider the 'backward shift operator' $S^*$ on $PW_a$ defined as $S^*f = \frac{f-f(0)}{z}$.
This is a bounded operator with a dense image.
We will denote the image of $S^*$ by $PW_a'$. Functions from $PW_a'$ appear naturally in relation with \Sch operators.

\begin{theorem}[\cite{MP3}]\label{tHBS}
An entire function $E = A+iB$ of Hermite-Biehler  class  corresponding to a \Sch equation from $\SS^2$ satisfies
\begin{equation}
A = \sin z + f \textrm{ and } B = \cos z + g, \label{eAB}
\end{equation}
where $f$ is an odd function from $PW_1'$ and $g$ is an even function from $PW_1'$ such that $f = S^*F, g = S^*G$ for some $F, G\in PW_1,\ F(0) = G(0)$.

Moreover, for any such $f, g\in PW_1'$ with $F(0) = G(0)$
there exists $\e > 0$ such that for any $c_1 \in \R,\ |c_1| < \e$ and $c_2 \in \R,\ |c_2| < \e$
the function $E = A+iB$ with
\begin{equation}
A = \sin z +c_1 f \textrm{ and } B = \cos z + c_2 g \label{eAB2}
\end{equation}
is an Hermite-Biehler function corresponding to a \Sch equation from $\SS^2$.
\end{theorem}

\begin{remark}\normalfont Note that, as follows from the proof in \cite{MP3}, the constant $\e$ can always be chosen as $1/\max(||f||_2, ||g||_2)$.

The statement of the theorem implies the asymptotics of the spectra \eqref{e000a}, \eqref{e000} and of the pointmasses of spectral measures \eqref{asymp1}, \eqref{asymp2} and is in fact equivalent to those asymptotics via Parseval's theorem
and the equivalence of the corresponding de Branges chains to PW chains.
\end{remark}

\section{Completeness, Gap and Type Problems}\label{GapType}

\subsection{Beurling-Malliavin densities and the radius of completeness}\label{BM}

If $\{I_n\}$ is a sequence of disjoint intervals on $\R$, we call it short if
$$\sum \frac{|I_n|^2}{1+\dist^2(0,I_n)} < \infty$$
and long otherwise.

If $\L$ is a sequence of real points define its exterior Beurling-Malliavin (BM) density (effective BM density) as
$$D^*(\L)=\sup\{ d\ |\ \exists\textrm{ long  }\{I_n\}\textrm{ such that }\#(\L\cap I_n)\geq d|I_n|\},\ \forall n\}.$$

For a non-real sequence its density can be defined as $D^*(\L)=D^*(\L')$ where $\L'$ is a real sequence $\lan'=\frac 1{\Re \frac 1{\lan}}$,
if $\L$ has no imaginary points, or as $D^*(\L)=D^*((\L+c)')$, with a properly chosen real constant $c$, otherwise.

A dual definition is used to introduce the interior BM density:
$$D_*(\L)=\inf\{ d\ |\ \exists\textrm{ long  }\{I_n\}\textrm{ such that }\#(\L\cap I_n)\leq d|I_n|\},\ \forall n\}.$$

Both densities play important role in Harmonic Analysis by appearing in a number of fundamental results. 
Their applications were recently extended into the area of spectral problems for differential operators via the methods discussed in this paper.
As an example, let us recall the original appearance of $D^*$ in the solution of a completeness problem by Beurling and Malliavin.

For any complex sequence $\L$ its radius of completeness is defined as
$$R(\L) =  \sup\{ a\ |\ \EE_\L=\{e^{i\l z}\}_{\l\in\L} \textrm{ is complete in }L^2(0,a)\}.$$

One of the fundamental results of Harmonic Analysis is the following theorem (see \cite{CBMS} for history and further references).

\begin{theorem}[Beurling and Malliavin, around 1961, \cite{BM1, BM2}]\label{bigBM}
Let $\L$ be a discrete sequence. Then
$$R(\L)=2\pi D^*(\L).$$
\end{theorem}

We will return to this result and the exterior density when we discuss Horvath' theorem in Section \ref{sHorvath}.
The other density, $D_*(\L)$, which has recently made a new appearance in the area of the Gap and Type Problems (see \cite{Polya, Gap, Type, CBMS}) are used in our statements below.
Note that for subsequences of sequences close to arithmetic progressions, the two densities are related to each other in a rather simple way.
In particular for sequences $\L$ satisfying \eqref{asymp1} or \eqref{asymp2}, for any $\Phi \subset \L$,
$$D^*(\Phi)+D_*(\L\setminus\Phi)=D^*(\L)=D_*(\L) =1.$$

\subsection{Spectral gaps, types and sign changes}\label{Gap}

We will call a lower semi-continuous function $W:\R\to [1,\infty]$ a weight on $\R$. We will say that a measure $\mu$ on $\R$ is $W$-finite if
$$ \| \mu \|_W=\int Wd|\mu|<\infty.$$
Note that $W$-finite measures are forced to be supported on the subset of $\R$ where $W$ is finite.

If $W$ is a weight we define its type $T_W$ as
$$T_W=\sup \{ a|\textrm{ $\exists$ $W$-finite non-zero measure $\mu$ with a spectral gap $[-a,a]$}\}.$$
This is one of several equivalent definitions for $T_W$, see \cite{CBMS}. The type of $W$ is used in applications such as problems on completeness of exponential systems in $L^p$ and Bernstein's spaces.

If $\nu$ is a real measure on $\R$ we will denote by $\nu^\pm$ its positive and negative parts.
We will say that $\nu$ has a spectral gap $(-a,a)$ if $\int fd\nu=0$ for all $f\in PW_a\cap L^1(|\nu|)$.
Note that for finite measures this property coincides with $\hat{\nu}=0$ on $(-a,a)$.

If $A$ and $B$ are two closed subsets of $\R$, let $\MM^W_a(A,B)$ be the class of all non-zero $W$-finite real measures $\sigma$ with a spectral gap $[-a,a]$ such that $\supp \sigma^+\subset A$ and $\supp \sigma^-\subset B$.

The following statement is a combination of Lemmas 13 and 16 from \cite{Det}.

\begin{lemma}\label{MiP}
If $\MM^W_a(A,B)$ is non-empty it contains a discrete measure $\nu$, whose positive and negative parts have interlacing supports.
\end{lemma}

The following version of the Type theorem is Theorem 36 from \cite{CBMS}
\begin{theorem}\label{Bmain}
$$T_W=\pi\sup \left\{d\ : \sum\frac{\log W(\lan)}{1+\lan^2}<\infty\textrm{ for some $d$-uniform sequence }\L \right\},$$
if the set is non-empty, and $0$ otherwise.
\end{theorem}

The notion of $d$-uniform sequences first appeared in \cite{Gap}.
For our results, we do not need a full definition of a $d$-uniform sequence, since the sequences we are concerned with are separated.
For a separated sequence $\L$, i.e., a discrete sequence such that 
$|\lan-\l_{n-1}|>c>0$ for all $n$, $\L$ is $d$-uniform iff $D_*(\L)=d$, see
for instance Example 1 on page 27 of \cite{CBMS}. This implies

\begin{corollary}[\cite{MP3}]\label{type}
$$T_W\geq \pi\sup \left\{d\ : \sum\frac{\log W(\lan)}{1+\lan^2}<\infty\textrm{ for some separated }\L,  \  D_*(\L)=d\right\}.$$
\end{corollary}

From one of the results of \cite{Det} (or Theorem 17, Chapter 4, in \cite{CBMS}) we deduce
\begin{corollary} [\cite{MP3}] \label{forTeven} If $W$ is a weight, $A$ is a separated sequence and $B$ is any closed set then
$$\sup\{a|\ \MM^W_a(A,B)\neq \emptyset\}\leq 2\pi D_*(A).$$
\end{corollary}

We denote $\log_-x=\max(0,-\log x).$
We will also need the following
\begin{corollary} [\cite{MP3}] \label{l2} Let $X=\{x_n\}$ be a separated sequence of real numbers.
Let $\tau=\sum c_n\delta_{x_n}$ be a discrete measure with spectral gap $(-a,a)$. Then for any $d<a$ there exists $\Phi\subset X$, $\pi D_*(\Phi)>d$, such that
$$\sum_{x_n\in\Phi} \frac{\log_- c_n}{1+n^2}<\infty.$$
\end{corollary}

\subsection{Horvath' theorem}\label{sHorvath}
The following theorem is one of the main results of \cite{Horvath}.
Here we formulate it
in an equivalent form using our version of the $m$-functions (after the square root transform).

\begin{theorem}[\cite{Horvath}] \label{Horvath}
Let $\L=\{\lan\}$ be a sequence of distinct non-zero complex numbers, $0 \leq a \leq 1$.
The following statements are equivalent:
\begin{enumerate}
\item For any \Sch operator $L\in \SS^2$, if $\ti L\in \SS^2$ is such that $q=\ti q$ on $(0,a\pi)$ and $m_-=\ti m_-$ on $\L$, then $L$ and $\ti L$ coincide identically.
\item The system of exponentials $\left\{ e^{i\gamma z}| \gamma \in \L\cup\{\ast,\ast\} \right\}$ is complete in \newline $L^2(0,(2-2a)\pi)$.
\end{enumerate}
\end{theorem}

The notation $\{*,*\}$ in the statement stands for any two points in $\C\setminus\L$.
A version of the above theorem is proved in \cite{Horvath} for all $1 \leq p \leq \infty$.
Here we treat only the case $p=2$, although a similar argument can be applied to other $p$.
A simple proof for the "un-mixed" case $a=0$ is given in \cite{BBP}.
Alternatively, this theorem can be deduced from Theorem \ref{tHBS} and Lemma \ref{l5}, details can be found in \cite{MP3}.

Horvath's theorem establishes equivalence between mixed spectral problems for \Sch operators and the Beurling-Malliavin problem on completeness of exponentials in $L^2$ spaces discussed in Section \ref{BM}.
Combining Theorems \ref{Horvath} and \ref{bigBM} we obtain the following statement, which is the sharpest possible result formulated in terms of the density of the defining sequence $\L$.

\begin{corollary}[\cite{MP3}]
Let $\L=\{\lan\}$ be a sequence of distinct non-zero complex numbers, $0 \leq a \leq 1$.
The following statements are equivalent:
\begin{enumerate}
\item Any two \Sch operators $L$ and $\ti L$, such that $q=\ti q$ on $(0,d\pi)$ for some $d > a$ and $m_-=\ti m_-$ on $\L$, coincide identically.
\item $D^*(\L) \geq 1-a$.
\end{enumerate}
\end{corollary}

\section{Applications to mixed spectral problems}\label{uncertainty}
Another result we want to include here is a 'parametrization' of counterexamples (indeterminate operators) in the so-called three-interval case of the mixed spectral problem.
The theorem by Hochstadt and Lieberman \cite{HL} says that knowing the potential on one-half of the interval $(0,\pi/2)$ and knowing one of the spectra allows one to recover the operator uniquely. The result is precise in the sense that the knowledge of the spectrum minus one point, or of the potential on $(0,\pi/2-\e)$, is insufficient.
Further results by del Rio, Gesztesy, Simon and Horvath focus on the case of two intervals, i.e., the case when the potential is known on $(0,a),\ 0<a<\pi$ and unknown on $(a,\pi)$. 
The logical problem is to try to replace
$(0,a)$ with other subsets of $(0,\pi)$, starting with a natural next step of two intervals $(0,a)\cup (b,\pi),\ 0<a<b<\pi$. However, such attempts immediately
meet the following elegant counterexample by Gesztezy and Simon in \cite{SG2}. 

\begin{counterexample} [\cite{SG2}]
Let the potential $q$ of a \Sch operator $L$ satisfy $q(x)=q(\pi-x)$ for all $x\in(0,\frac \pi 2-\e)$. Let $\ti L$ be the operator with the reflected potential $\ti q(x)=q(\pi-x)$ for all $x$. It is well known that
then the spectra of the two operators, with any pair of symmetric boundary conditions at the endpoints, will coincide. Nonetheless, the operators are not identical, unless $q(x)=q(\pi-x)$ for $x\in(\frac \pi 2-\e, \frac \pi 2)$.
Thus, having almost complete direct information (knowing the potential on the two intervals of total length $\pi-2\e$) and a spectrum is insufficient to recover the operator.
\end{counterexample}

This counterexample shows that the two-interval problem is substantially different from the three-interval case.
Until recently this counterexample for the three-interval problem was the only one existing in the literature, prompting some of
the experts to ask if it was, in some sense, the only counterexample of that kind. Theorem \ref{3int} in Section \ref{s3int} is a result describing all possible counterexamples
for the three-interval problem in terms of the spectral measures of the operators. It follows from our 'parametrization' that the pairs of operators with
'almost symmetric' potentials provide only a  submanifold of all pairs of indeterminate operators.
All other such pairs, however, have to be close to symmetric in terms of the asymptotics of the spectral measure, see Section \ref{s3int}.
In the opposite direction, we show that
there exists a large class of operators which are uniquely determined by their potentials on two intervals and a proper part of the spectrum.
This raises the natural question of describing the operators with such uniqueness properties discussed in Section \ref{P2}.

\subsection{Fourier gaps and Schr\"odinger operators}\label{FSGap}
We formulate results of this section for a slightly broader class $\SS^1$ of \Sch operators with summable potentials.

\begin{lemma}[\cite{MP3}]\label{l5}
Let $L,\ti L\in \SS^1$. Then  $q=\ti q$ on $(0,a\pi)$ iff $\mu_--\ti\mu_-$ has a spectral gap  $(-2a,2a)$. Similarly,
$q=\ti q$ on $(b\pi,\pi)$ iff $\mu_+-\ti\mu_+$ has a spectral gap  $(-2(1-b),2(1-b))$.
\end{lemma}

\begin{remark}\normalfont It is a well known statement in the area of the gap problem that a real measure $\mu$ has a spectral gap $(-a,a)$
if and only if $H\mu(iy)=O(e^{-ay})\textrm{ as }y\to\infty,$ see for instance Lemma 4.5 in \cite{Polya}. In this statement $O(e^{-ay})$ can be replaced with
$o(e^{-y(1-\e)a})$ and the positive imaginary half-axis with negative.
\end{remark}

\begin{lemma}[\cite{MP3}]\label{l6} Let $L,\ti L\in \SS^1$ and let $k,\ti k$ be the corresponding Krein functions. Then $q=\ti q$ on $(0,a\pi)$ iff the measure
$(k-\ti k) dx$ has a spectral gap  $(-2a,2a)$.
\end{lemma}

As was mentioned above, the presence of spectral gap can be equivalently reformulated in terms of decay along $i\R$:

\begin{corollary}[\cite{MP3}]\label{c01}Consider two operators $L$ and $\ti L$ from $\SS^2$ and let $k,\ \ti k$ be their Krein functions. Then
the potentials $q$ and $\ti q$ are equal on $(0,a\pi)$ iff
$$H(k-\ti k)(iy)=O(e^{-2ay})\textrm{ as }y\to\infty$$
iff
$$H(k-\ti k)(iy)=o(e^{-2(1-\e)ay})\textrm{ as }y\to\infty$$
for any $\e>0$.

\end{corollary}

\begin{remark}\normalfont It follows from the proof of Lemma \ref{l5} that  $q=\ti q$ on $(0,a\pi)$ iff
$\nu_- -\ti \nu_-$  has a spectral gap  $(-2a,2a)$, where $\nu_-$ is the spectral measure
corresponding to the Neumann boundary condition at 0.  The same fact follows from Lemma \ref{l6} because $\pi - k$ and $\pi - \ti k$ are
Krein functions for $\nu_-$ and $\ti \nu_-$.

More generally, any self-adjoint boundary condition at 0 or $\pi$ can be used in this statement with the same proof.
\end{remark}

Together with the above this produces the following statement, which is one of the main results of \cite{S}. It is formulated in \cite{S} in an equivalent form,
using a different $m$-function (before the square root transform).

\begin{corollary}[\cite{MP3}]\label{c1}
The potentials $q$ and $\ti q$ coincide on $(0,a\pi)$ iff the $m$ functions satisfy
$$m_+(iy)-\ti m_+(iy)=O(e^{-2ay})\textrm{ as }y\to\infty.$$
\end{corollary}

\subsection{The size of uncertainty}\label{2int}
Recall that for an operator $L$ we denote by $\sigma_{DD}=\{\lan\}$ and $\sigma_{ND}=\{\eta_n\}$ its spectra after the square root transform. 
We enumerate the sequences as described in Section \ref{spectra}.

Let $I=\{I_n\}_{n\in\Z}$ be a sequence of intervals on $\R$. 
Consider the set $\SS_I$ of those operators $L\in \SS^2$ for which
$\sigma_{DD}$ and $\sigma_{ND}$ lie in the union of $I_n$. 
Following our discussion in the introduction,
we can ask what part of the necessary spectral information we are given by this inclusion condition?
Let us define the size of uncertainty for the sequence of intervals $I=\{I_n\}$ as the number $U(I)$
equal to the infimum of $a$ such that knowing the potential of an operator $L\in\SS_I$ on $(0,a\pi)$
one can recover $L$ uniquely. Theorem \ref{t0} below gives the following formula for the size of uncertaity:
$$U(I)=\pi \sup \left\{D_*(\Phi)\ :\ \  \sum_{n\in\Phi}\frac{\log_-|I_n|}{1+n^2}<\infty\right\}.$$

\no Let us now make our statements more precise.

\begin{theorem}[\cite{MP3}]\label{t0}
Let $\{\e_n\}_{n\in\N}$ be positive numbers,  $a\in [0,1)$. The following statements are equivalent:

\begin{enumerate}
\item Any two Schr\"odinger operators $L$ and $\ti L$ satisfying
\begin{equation}|\lan-\ti\lan|<\e_{2n},\ |\eta_n-\ti\eta_n|<\e_{2n+1}, \ n\in\N,\label{e1}\end{equation}
 and $q(x)=\ti q(x)$ on $(0,d\pi), d>a$ must coincide identically, i.e.,  $q=\ti q$ a. e. on $(0,\pi)$.
\item Any sequence of distinct integers $\Phi\subset \Z$ such that
 \begin{equation}\sum_{n\in\Phi\cap\N}\frac{\log_-\e_n}{1+n^2}<\infty,\label{e2}\end{equation}
satisfies $ D_*(\Phi)\leq a$.
\end{enumerate}
\end{theorem}

In particular we obtain the following 'uncertainty' version of Borg's theorem
\begin{corollary}[\cite{MP3}]
Let $\{I_n\}$ be a sequence of intervals on $\R$, $U=\cup I_n$. The following statements are equivalent:

\begin{enumerate}
\item The condition \begin{equation}\sigma_{DD}\cup \sigma_{ND}\subset U\label{eUP}\end{equation}
together with the values of the potential $q$ on $(0,\e)$ for any $\e>0$ determines an operator $L\in \SS^2$
uniquely.
\item For every $\Phi\subset \N$ satisfying
$$\sum_{n\in\Phi}\frac{\log_-|I_n|}{1+n^2}<\infty,$$
there exists a long sequence of intervals $\{J_n\}$ in $\R_+$ such that
$$\frac{\#(\Phi\cap J_n)}{|J_n|}\rightarrow 0$$
as $n\to\infty$.
\end{enumerate}
\end{corollary}

The first statement can be equivalently formulated as follows: 
if $L,\ti L\in \SS^2$ satisfy \eqref{eUP} and $q = \ti q$ on $(0,\e)$ for some $\e>0$ then $L\equiv \ti L$.
Recall that the definition of long sequences of intervals (in the sense of Beurling and Malliavin) was given in Section \ref{BM}.

\subsection{Three-interval statements}\label{s3int}

Our goal in this section is to describe all possible counterexamples in the three interval problem.
To formulate Theorem \ref{3int} below we will need some preparation.

Let $\mu$ be a positive measure on $\R$ such that $PW_a\subset L^2(\mu)$ for all $a>0$. It is well known, see for instance \cite{OS},
that any spectral measure of a \Sch operator from $\SS^2$ satisfies this condition. We define the type of $\mu$ as
$$T_\mu=\inf \{a|\ PW_a\textrm{ is dense in }L^2(\mu)\}.$$

Let now $\L$ be a sequence satisfying the asymptotics \eqref{e000a}. Denote by $\eta$ the counting measure of $\L$.
It follows from the Type Theorem of \cite{Type, CBMS}, and in fact from much earlier results on the type problem, that $T_\eta=1$.
Moreover, as follows for instance from the results of  \cite{OS}, $L^2(\eta)=PW_1$ (as sets).

Let now $0<c,d<1$ be real constants. In our next statement we will need  two even functions $f,g$,
$$f\in L^2(\eta)\ominus PW_{c},\ \  g\in L^2(\eta)\ominus PW_{d}$$
which satisfy certain asymptotics, see \eqref{e3}. As we will discuss after the statement, such functions form
dense subsets in the corresponding infinite-dimensional subspaces of $L^2(\eta)$. Recall that by $\Gamma(\L)$ we denote
the characteristic sequence of $\L$, or equivalently the sequence of pointmasses of the even operator, see Section \ref{even}.

\begin{theorem}[\cite{MP3}]\label{3int} Let $\L$ be a sequence satisfying the asymptotics \eqref{e000a}, $\Gamma(\L)=\{\gamma_n\}$. Denote by $\eta$ the counting measure of $\L$. Let
 $a,b\in [0,1/2)$ be arbitrary constants

\no 1) for any two real even functions $f,g$,
$$f\in L^2(\eta)\ominus PW_{2a},\ \  g\in L^2(\eta)\ominus PW_{2b}$$
 satisfying
\begin{equation}|n|f(\lan)=a_n,\ |n|g(\lan)=a_n\left(1+\frac{b_n}{|n|+1}\right),\label{e3}\end{equation}
 for some $ a_n,b_n\in l^2$, $b_n\geq -|n|-1$,
the measures $\mu=\sum \alpha_n\delta_{\lan}$ and  $\ti\mu=\sum \ti\alpha_n\delta_{\lan}$, where
$$\alpha_n=\frac {f(\lan)}2 + \sqrt{\frac{f^2(\lan)}4 + \gamma^2_n\frac {g(\lan)}{f(\lan)}},$$
\begin{equation}\ti\alpha_n=-\frac {f(\lan)}2 + \sqrt{\frac{f^2(\lan)}4 + \gamma^2_n\frac {g(\lan)}{f(\lan)}}\label{e5}\end{equation}
when $a_n\neq 0$ and $\alpha_n=\ti\alpha_n=1+c_n/n$ for some $c_n\in l^2$ when $a_n=0$,
are spectral measures for two \Sch  operators $L,\ti L\in \SS^2$  such that
$$\sigma_{DD}=\ti\sigma_{DD}=\L,\ q=\ti q\textrm{ on }[0,a\pi]\cup [(1-b)\pi,\pi].$$

\no 2) For any two spectral measures of Schr\"odinger operators $L,\ti L$ on $[0,\pi]$ such that
$$\sigma_{DD}=\ti\sigma_{DD}=\L,\ q=\ti q\textrm{ on }[0,a\pi]\cup [(1-b)\pi,\pi],$$
their pointmasses must satisfy
$$\alpha_n-\ti\alpha_n=f(\lan), \ \ti\alpha^{-1}_n-\alpha^{-1}_n=\gamma_n^{-2}g(\lan),$$
for some even functions $$f\in L^2(\eta)\ominus PW_{2a},\ g\in L^2(\eta)\ominus PW_{2b},$$
satisfying \eqref{e3} with some $ a_n,b_n\in l^2$, $b_n\geq -|n|-1$.

\end{theorem}

\begin{remark}\normalfont
 In regard to the functions $f$ and $g$ from the statement, note that since $2a, 2b <1$ and $L^2(\eta)=PW_1$, the orthogonal complements of $PW_{2a}$
 and $PW_{2b}$ in $L^2(\eta)$ are large infinite-dimensional subspaces.

\no If a function $h$ belongs to $ L^2(\eta)\ominus PW_{2a}$, then $p(z)=\textrm{Re }h(z)-\textrm{Re }h(-z)$ is an odd real function from the same subspace.
The function $f(z)=p(z)/z$ will satisfy  the conditions for $f$. Moreover, the set of such functions is dense in the space of
all even  functions from $L^2(\eta)\ominus PW_{2a}$.

\no To choose $g$, assume for instance that $a\geq b$. Choose two more real odd functions $p_1, p_2$ in $ L^2(\eta)\ominus PW_{2a}$ as described above. Then the function
$r(z)=p_1(z)p_2(0)-p_1(0)p_2(z)$ will have a double zero at 0. The function $ g=f+r/z^2$ will satisfy the conditions for $g$. If $a< b$ one
may proceed in the opposite direction, first choosing $g$ and then $f$.

\end{remark}

\no Theorem \ref{3int} describes all possible counterexamples for the three-interval problem. As we can see, this set is richer  than previously thought, with the original counterexample from \cite{SG2} discussed in the introduction corresponding to the case $f=g$, i.e., $b_n\equiv 0$ in the notations of the statement. On the other hand, all of the counterexamples must be close to the original in the sense that the difference between
$f$ and $g$ has higher order of decay in comparison with either function.

\subsection{Uniqueness in the 3-interval problem near the even operator}\label{3intU}

\no In this section we will look at the 3-interval problem from a slightly different point of view. We will ask if there exist
operators $L$ uniquely determined by the '3-interval' information. Our result here can be viewed as a set-up for the problem
of description of such operators discussed in the next section.

\no Let $\mu=\mu_-=\sum \alpha_n\delta_{\lan}$ denote the spectral  measure of an  operator $L$ from $\SS^2$ as defined in Section \ref{functions}.
For $\L\subset \Z$ denote
$$D'_*(\L)=\max\{D_*(\Phi)|\ \Phi\subset\L, \{n,n+1\}\not\subset \Phi, \forall n\in \Z\}.$$
Once again, for a discrete sequence $\L\subset\R$ we denote by $\Gamma(\L)$ its characteristic sequence
from Section \ref{even}.

\begin{theorem}\label{tUoper} Let $ 0\leq a<1/2$ and let $\{\e_n\}_{n\in\N}$ be a sequence of positive numbers such that for
$\Phi\subset\N$ the inequality
\begin{equation}\sum_{n\in\Phi}\frac {\log_-\e_n}{1+n^2}<\infty\label{e22}\end{equation}
implies
$$D'_*(\Phi)\leq a.$$
Then
any operator $L\in \SS^2$, satisfying
\begin{equation}|\alpha_n-\gamma_n|<\e_n,\label{e21}\end{equation}
where $\{\gamma_n\}=\Gamma(\sigma_{DD})$, is uniquely determined by its spectrum $\sigma_{DD}$ and its potential
$q$ on $(0,d\pi)\cup ((1-d)\pi,\pi)$ for any $d>a$. I.e., if $L$ satisfies \eqref{e21} then any $\ti L\in \SS^2$ with $\ti\sigma_{DD}=\sigma_{DD}$ and
$\ti q=q$ on $(0,d\pi)\cup ((1-d)\pi,\pi)$ for some $d>a$ must be identical to $L$.

\end{theorem}

\begin{remark}\normalfont The last statement displays an interesting phenomenon. As we can see from the original counterexample, if an operator is close to even (but not even)
in terms of potential (in the sense of direct problem), namely if its potential is even on $[0,d]\cup [1-d,\pi]$, then it is not uniquely defined by the corresponding
mixed spectral data. However, if an operator is close to even in terms of the spectral measure (in the sense of inverse problem) then it is uniquely defined by the mixed spectral data. This property is yet to be fully understood. One of its particular consequences is a statement typical for the
area of the Uncertainty Principle: an operator cannot be close to even (without being even) simultaneously with respect to its potential and spectrum, in the sense of the last statement.

In place of the condition $\ti\sigma_{DD}=\sigma_{DD}$, one may instead assume that the intersection of the spectra has sufficiently large density. This alternative condition follows from the proof of the theorem in \cite{MP3}.
\end{remark}

\section{Examples and open problems} \label{examples}
In this section we provide explicit examples for spectral problems using methods and results we stated throughout the paper.

We start with examples of spectral problems.

\subsection{Constant Hamiltonians}

With Theorems \ref{t3} and \ref{t007} or Theorems \ref{tPer} and \ref{tPer2}, we obtain solutions to spectral problems, inverse and direct, for systems
with constant Hamiltonians.

\begin{example}[\cite{etudes}]
Let the constants $h_1 > 0, h_2  >0 , g \in \R$ satisfy   $h_1 h_2 = 1 + g^2$, and
$$\HH = \begin{pmatrix} h_1 & g \\ g & h_2 \end{pmatrix} \in {\rm SL}(2,\mathbb R).$$
Then $\mu = h_1^{-1}$.
\end{example}

\subsection{Even $2\pi-$periodic spectral measures}
Using Theorem \ref{tPer} or \ref{tONP}, we solve the inverse spectral problem for this particular spectral measure $d\mu = (1+\cos x) dx$.
\begin{example}[\cite{etudes}]
For the spectral measure $d\mu = (1+\cos x) dx$, $h$ is a step function of step size $\frac{1}{2}$, and the steps take value \begin{equation*}
    1, \frac{1}{3}, \frac{2}{3}, \frac{2}{5}, \frac{3}{5}, \ldots, \frac{n}{2n - 1}, \frac{n}{2n + 1}, \frac{n + 1}{2n + 1}, \ldots
\end{equation*}
on intervals 
\begin{equation*}
    \left[0, \frac{1}{2}\right), \left[ \frac{1}{2}, 1 \right), \left[ 1, \frac{2}{3} \right), \ldots
\end{equation*}

or equivalently,
$$h^\mu_n = \frac{ \left( (-1)^n(2n + 3) + 1 \right)^2}{8(n + 1)(n + 2)} \quad \text{on} \quad \left[ \frac{n}{2}, \frac{n + 1}{2} \right).$$
\end{example}

Similarly, we can solve the ISP for the spectral measure $d\mu = (1-\cos x) dx$.
\begin{example}[\cite{etudes, Direct}]
For the spectral measure $d\mu = (1 - \cos x) dx$, $h$ is a step function of step size $\frac{1}{2}$, and and satisfies
$$h^\mu_n = \frac{(n + 1)(n + 2)}{2} \quad \text{on} \quad \left[ \frac{n}{2}, \frac{n + 1}{2} \right).$$
\end{example}

We obtained the following example by using Theorem \ref{directalg}.
\begin{example}[\cite{Direct}]\label{GeometricGrowth}
Consider the diagonal det-normalized Hamiltonian satisfying $h^\mu_n = a^n$ on $[\frac{n}{2}, \frac{n + 1}{2})$ where $a > 0$.

This spectral measure corresponds to the Geronimus polynomials: let $\alpha \equiv \frac{1 - a}{1 + a}$, on the unit circle $\mu$ has an absolutely continuous part \begin{equation*}
    w(x) = \frac{1}{|1 + \alpha|} \frac{\sqrt{1 - \alpha^2 - \cos^2 \frac{x}{2}}}{\sin \frac{x}{2}},
\end{equation*}
supported on $\left[ 2\arcsin{|\alpha|}, 2\pi - 2\arcsin{|\alpha|} \right]$, and a singular part \begin{equation*}
    \frac{2}{|1 + \alpha|^2} \left( |\alpha + \frac{1}{2}|^2 - \frac{1}{4} \right) \delta_{z = 1},
\end{equation*}
if $\alpha > 0$, or equivalently $a < 1$.
\end{example}

\subsection{General $2\pi-$periodic measures}

Using Theorems \ref{tPer} and \ref{tPer2} or \ref{tONP} and \ref{tONP2}, we can solve the inverse spectral problem for the spectral measure $d\mu = (1 \pm \sin x) dx$.

\begin{example}[\cite{etudes}]
For the periodic spectral measure $d\mu = (1 + \sin x) dx$, we can use Theorem \ref{tPer} or \ref{tONP} to find $h^\mu$, Theorem \ref{tPer2} or \ref{tONP2} to find $g^\mu$. The Hamiltonian takes value \begin{equation*}
\begin{aligned}
&\begin{pmatrix}
1 & 0 \\
0 & 1
\end{pmatrix}, \quad
\begin{pmatrix}
5/3 & -4/3 \\
-4/3 & 5/3
\end{pmatrix}, \quad
\begin{pmatrix}
4/3 & -5/3 \\
-5/3 & 17/6
\end{pmatrix}, \quad
\begin{pmatrix}
4/5 & -1 \\
-1 & 5/2
\end{pmatrix}, \\
&\begin{pmatrix}
13/15 & -2/3 \\
-2/3 & 5/3
\end{pmatrix}, \quad
\begin{pmatrix}
25/21 & -22/21 \\
-22/21 & 37/21
\end{pmatrix}, \quad
\begin{pmatrix}
8/7 & -9/7 \\
-9/7 & 65/28
\end{pmatrix}, \ldots
\end{aligned}
\end{equation*}
on intervals 
\begin{equation*}
    \left[0, \frac{1}{2}\right), \left[ \frac{1}{2}, 1 \right), \left[ 1, \frac{2}{3} \right), \ldots
\end{equation*}
\end{example}

\begin{example}[\cite{etudes}]
For the spectral measure $d\mu = (1 - \sin x) dx$, we can again obtain the Hamiltonians similar to the previous example: 
\begin{equation*}
\begin{aligned}
&\begin{pmatrix}
1 & 0 \\
0 & 1
\end{pmatrix}, \quad
\begin{pmatrix}
5/3 & 4/3 \\
4/3 & 5/3
\end{pmatrix}, \quad
\begin{pmatrix}
4/3 & 5/3 \\
5/3 & 17/6
\end{pmatrix}, \quad
\begin{pmatrix}
4/5 & 1 \\
1 & 5/2
\end{pmatrix}, \\
&\begin{pmatrix}
13/15 & 2/3 \\
2/3 & 5/3
\end{pmatrix}, \quad
\begin{pmatrix}
25/21 & 22/21 \\
22/21 & 37/21
\end{pmatrix}, \quad
\begin{pmatrix}
8/7 & 9/7 \\
9/7 & 65/28
\end{pmatrix}, \ldots
\end{aligned}
\end{equation*}
on intervals 
\begin{equation*}
    \left[0, \frac{1}{2}\right), \left[ \frac{1}{2}, 1 \right), \left[ 1, \frac{2}{3} \right), \ldots
\end{equation*}
\end{example}

\subsection{Non-periodic measures and periodization approach}

Using Theorem \ref{tTT}, we can explicitly solve the ISP for $\mu$ being the Lebesgue measure plus a point mass.

\begin{example}[\cite{etudes}]\label{pointmassexample}
\label{pointmass}
Consider the spectral measure $d\mu = \alpha + \beta\pi\delta_0,\ \a > 0,\b \geq 0$.
Using Theorem \ref{tTT} we find
$$\psi_t = \frac{\sqrt{2\pi}}{2\pi\alpha+2\pi t\beta} 1_{(-t,t)},$$
and
$$k_t(0)=\frac1{\sqrt{2\pi}}\int\psi_t=\frac{2tc(t)}{\sqrt{2\pi}}=\frac1\pi\frac t{\alpha+t\beta}.$$
It follows that
$$h^\mu(t)=\pi\frac d{dt} k_t(0)=\frac d{dt}\frac t{\alpha+t\beta}=\frac\alpha{(\alpha+t\beta)^2}$$
and that  $$H(t)=\begin{pmatrix}  \frac\alpha{(\alpha+t\beta)^2}&0\\0&\frac{(\alpha+t\beta)^2} \alpha     \end{pmatrix}$$
is a unique diagonal Hamiltonian with spectral measure $\mu$.

\begin{remark}
Note that in the limiting case $\alpha\to 0$ (and $\beta=1$):
$$\mu(x)\to \pi\delta(x),\qquad h^\mu(t)\to \delta(t).$$
Indeed,
$$\int_0^\infty\frac\alpha{(\alpha+t)^2}~dt=-\frac\alpha{\alpha+t}\Big|_{t=0}^\infty=1$$
This shows how using our methods developed for PW measures one can solve the inverse problem for a non-PW measure $\mu = \pi\delta$, which is consistent with some convergence results mentioned in Section \ref{periodic} from \cite{PZ}.
\end{remark}
\end{example}

\begin{example}[\cite{etudes}]
Here is a generalization of the last example studied with different methods by H. Winkler in \cite{W}.

Let $\mu$ be a PW sampling measure. For $r\geq 0$, let $\mu_r=\mu+r\pi\delta$. 
Then $$h^{\mu_r}(t)=\frac{h^\mu(t)}{\left(1+r\int_0^th^\mu\right)^2}.$$

This yields a lot of examples for the ISP. Let $\mu = \alpha m$, $h^\mu = \alpha^{-1}$, then we have
$$h^{\mu_r}(t)=\frac{\alpha^{-1}}{\left(1+\alpha^{-1}rt\right)^2}=\frac{\alpha}{(\alpha+rt)^2}.$$

Note that $\mu_{x+y}$ can be obtained as $\mu_{x+y}=\mu_x+y\d_0$. For $h_x=h^{\mu_x}$ this result yields the equation
$$h_{x+y}(t)=\frac{h_x(t)}{\left(1+y\int_0^th_x\right)^2}$$ with the initial condition $h_0=h^\mu$.
\end{example}

\begin{example}[\cite{etudes}]
Consider the spectral measure $\mu=\alpha+\beta\pi\delta_\lambda$. Here $\alpha>0$, $\beta>0$, and $\lambda\in \R$ are parameters. 
Unlike Example \ref{pointmass}, if $\l\neq 0$, then the measure is not even. 

\begin{enumerate}[label=(\roman*)]
\item We want to find
$$h^\mu=\sqrt{\frac\pi2}\frac d{dt}\int\psi_t, $$
where
$$\hat\mu\ast \psi_t=1\qquad{\rm on}\quad (-t,t).$$
We can solve this convolution equation and obtain
$$h^\mu(t)=\frac d{dt}\left[\frac t\alpha-\frac\beta\alpha\frac1{\alpha+\beta t}\left(\frac{\sin \lambda t}\lambda\right)^2\right].$$
E.g., if $\alpha=\beta=1$, then
$$h^\mu(t)=\sin^2\lambda t+\left[\cos \lambda t-\frac{\sin \lambda t}{\lambda(1+t)}\right]^2.$$

\item We now turn to the computation of $g^\mu$ for
$$\mu=\alpha+\pi\beta\delta_\lambda$$
To compute $g^\mu$ we need
$$l_t(0)=\frac 1{2\pi i}\int_{-t}^t\psi_t*(\sigma \hat\mu).$$

\noindent Via direct calculations, we can solve $l_t$ and obtain
\begin{equation*}
\begin{aligned}
g^\mu(t)
&= -\frac{\beta}{\alpha \pi \lambda}
\Bigg[
1
- \frac{\alpha}{\alpha + \beta t}\cos(2\lambda t) \\
&\qquad
- \frac{\alpha\beta + 2\beta^2 t}{(\alpha + \beta t)^2}
    \frac{\sin(2\lambda t)}{2\lambda}
+ \frac{\beta^2}{(\alpha + \beta t)^2}
    \frac{\sin^2(\lambda t)}{\lambda^2}
\Bigg].
\end{aligned}
\end{equation*}

\noindent Let us simplify the above formula in the special case $\alpha=\beta=\lambda=1$:
\begin{equation*}
\begin{aligned}
g^\mu(t)
&= -\frac{1}{\pi}
\Bigg[
1
- \frac{\cos(2t)}{1+t}
- \frac{1+2t}{(1+t)^2}\,\frac{\sin(2t)}{2}
+ \frac{\sin^2 t}{(1+t)^2}
\Bigg] \\
&= \frac{(2t+1)(\cos(2t)+\sin(2t)) + 2\cos(2t) - 1}
       {2\pi(1+t)^2}
- \frac{1}{\pi}.
\end{aligned}
\end{equation*}
\end{enumerate}
\end{example}

\begin{example}[\cite{etudes}]
Let us further generalize the previous example and solve the inverse problem for
$$\mu=\alpha +\sum_n \pi\beta_n \delta_{\lambda_n}.$$
For its Fourier transform
we have
$$\hat\mu=\sqrt{2\pi}\alpha\delta_0+\sum_n\frac{\pi\beta_n}{\sqrt{2\pi}}e_{-\lambda_n},
$$
where $e_{s}(t)=e^{ist}$.

As before, let $\psi=\psi_t$ be the Fourier transform of the reproducing kernel $k_t(\cdot)=K_t(0,\cdot)\in \BB_t$. 
Theorem \ref{tTT} then gives
$$\hat\mu*\psi=1\textrm{ on }(-t,t)$$
which can be rewritten as
\begin{equation}\sqrt{2\pi}\alpha\psi(x) + \sum_n c_n(t)\frac {\pi\beta_n}{\sqrt{2\pi}}e_{-\lambda_n}(x)=1\label{e01}\end{equation}
where
$$c_n(t)=\int_\R e_{\lambda_n}\psi =\sqrt{2\pi}\hat\psi(-\lambda_n).$$
Solving \eqref{e01} for $\psi$ we get
\begin{equation}
	\psi(x)=\frac 1{\sqrt{2\pi}\alpha}-\sum_n\frac{\beta_n}{2\alpha}c_n(t)e_{-\lambda_n}(x).\label{e02}\end{equation}

Direct calculations yield
\begin{equation}\left(\alpha +\beta_n t\right) c_n +\sum_{k\neq n} \ \si_t(\lambda_n-\lambda_k)\beta_kc_k=\frac 2{\sqrt{2\pi}}\ \si_t(\lambda_n),\label{e03}
\end{equation}
where $$\si_t(x)=\frac{\sin (tx)}{ x}.$$

Let us introduce the following two matrices,
$$B = \diag{\beta_1, \beta_2, \ldots, \beta_n} \quad \text{and} \quad S_t=\left(\si_t(\lambda_j-\lambda_k)\right)_{1\leq j,k\leq n},$$
and two vectors
$$
L_t=\sqrt{\frac 2\pi}\left(\begin{array}{c}
	\si_t(\lambda_1)\\
	\si_t(\lambda_2)\\
	\vdots \\
	\si_t(\lambda_n)
\end{array}\right),\ {\rm and}\  C(t)=\left(\begin{array}{c}
	c_1(t)\\
	c_2(t)\\
	\vdots \\
	c_n(t)
\end{array}\right).
$$
Our calculations above lead to the following
\begin{theorem}[\cite{etudes}]\label{t01}
	$$\left(\alpha+S_t B\right)^{-1} L_t=C(t).$$
\end{theorem}

It follows that $h^\mu$ can now be calculated via the following 
formula, which concludes the solution in the case when $\mu$ is even.

\begin{corollary}[\cite{etudes}]
$$	h^\mu(t)=\frac 1{\pi \alpha}-\frac 1{\sqrt{2\pi}\alpha}\frac d{dt}<B(\alpha+S_t B)^{-1} L_t,L_t>.$$
\end{corollary}

In the general case, the off-diagonal term of $\HH$, $g^\mu$, can be calculated using the formulas in Theorem \ref{tPer2}. 

\end{example}

\subsection{Periodization approach}
We first apply the periodization approach on Example \ref{pointmass}.
\begin{example}[\cite{etudes, PZ}]
Let us consider a special case of Example \ref{pointmass} $\mu = \frac{1}{\sqrt{2\pi}} m + \sqrt{2\pi} \delta_0$, and we know the corresponding $h^{\mu}(t)$ is
$$h^{\mu}(t) = \frac{\sqrt{2\pi}}{(1 + 2t)^2}.$$

This measure satisfies the assumptions of both theorems \ref{PW} and \ref{decay}. We periodize the measure with $T = \pi, 2\pi, 4\pi, 8\pi$, and compare the $h^{\mu}_T$'s with the actual $h^{\mu}$. In Figure \ref{pointmasspic}, the black step functions are the $h^{\mu}_T$'s obtained from the periodizations, and the orange curve in every figure is the actual $h^{\mu}(t)$.

\begin{figure}[ht]
\centering
\includegraphics[width=0.4\textwidth]{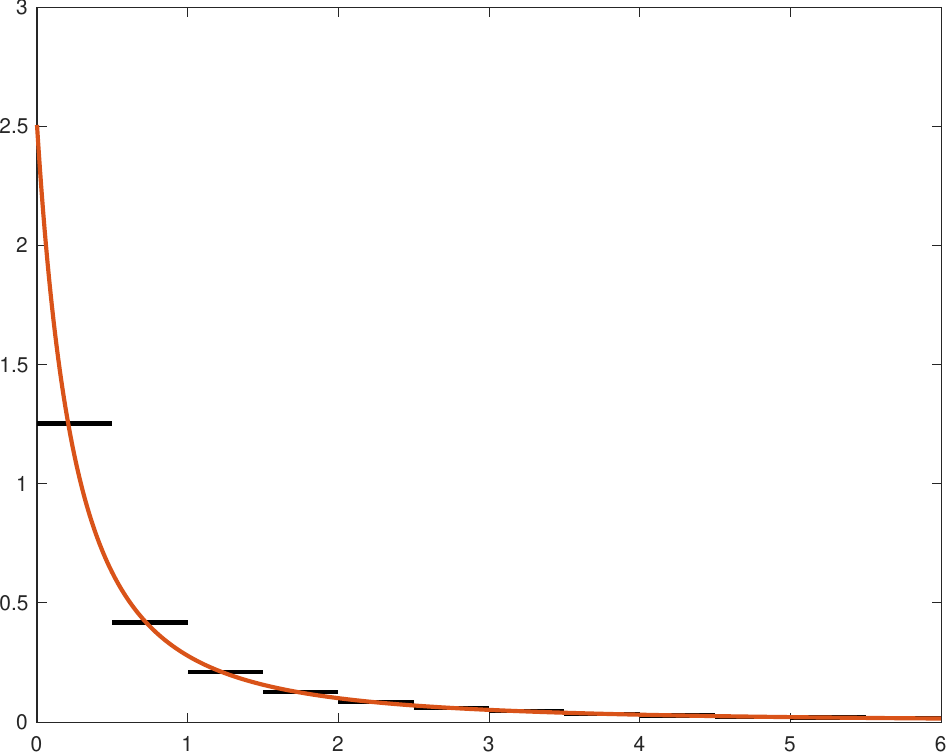}
\includegraphics[width=0.4\textwidth]{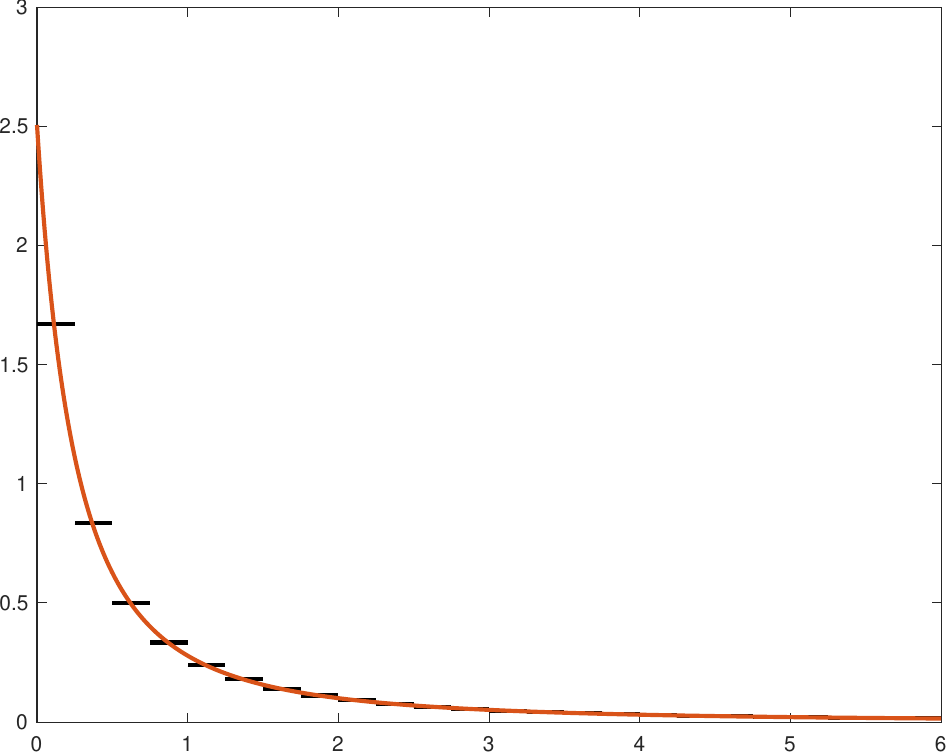}
\includegraphics[width=0.4\textwidth]{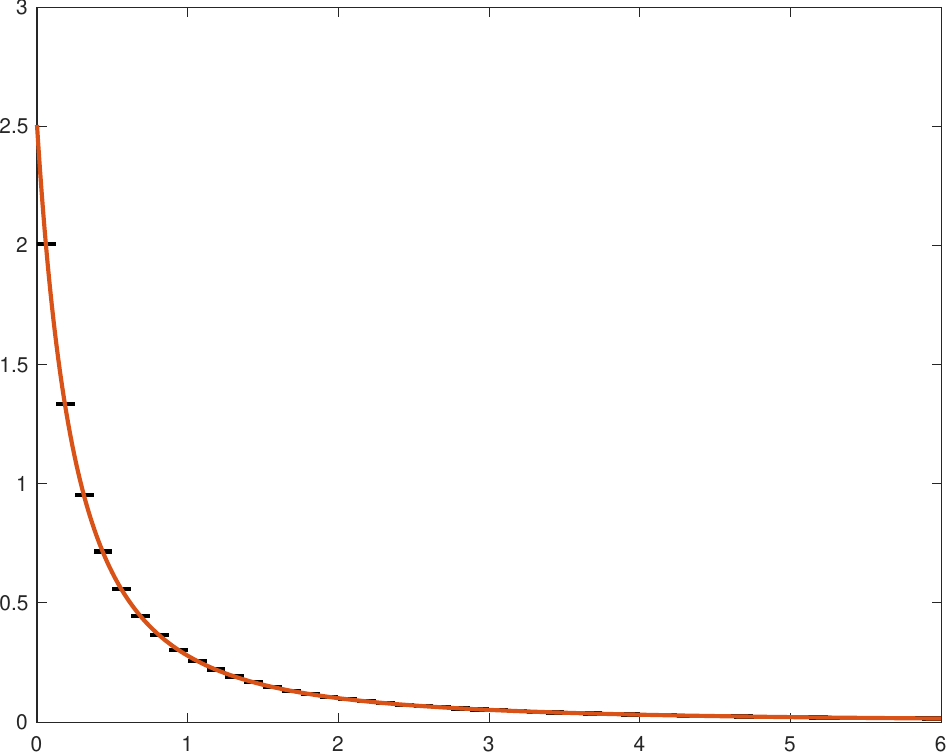}
\includegraphics[width=0.4\textwidth]{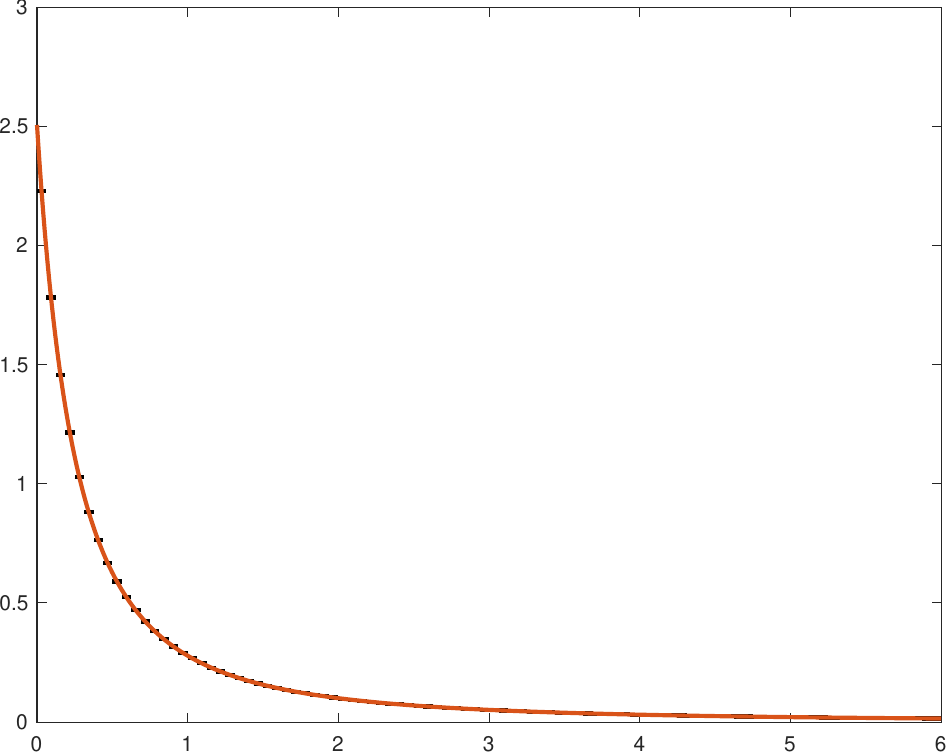}
\caption{From left to right: $h^{\mu}_{\pi}$, $h^{\mu}_{2\pi}$, $h^{\mu}_{4\pi}$, $h^{\mu}_{8\pi}$.}
\label{pointmasspic}
\end{figure}

For this example, we can obtain $h^{\mu}$ from $h^{\mu}_T$'s. The $2T-$periodization of $\mu$ is in the form $\mu_T = \sqrt{2\pi}\delta_0 + \frac{1}{\sqrt{2\pi}}m$ on $[-T, T]$e. Then the measure corresponding to $\mu_T$ on the unit circle $\mathbb{T}$ is $\mu_T^\mathbb{T} = \frac{\pi + T}{\sqrt{2\pi}} \left( \frac{\pi}{\pi + T} \delta_0 + \frac{T}{\pi + T} \frac{m}{2\pi} \right)$, where $\delta_0$ here stands for the unit point mass at $\theta = 0$ on $\mathbb{T}$. We know the orthonormal polynomials corresponding to this measure, for instance from \cite{Simon}. Therefore by theorem \ref{tONP}, the $n$-th step of $h^{\mu}_T$ corresponding to $\mu_T$ takes value \begin{equation*}
    \frac{\sqrt{2\pi}T^2}{(n\pi + T)(n\pi + T + \pi)}.
\end{equation*}

For a given $T$, in the places of $n$, plug in $\frac{2Tt}{\pi} + 1$. This produces the function whose graph passes through all the right endpoints of the steps for $h^{\mu}_T$. Denoting this function by $\Tilde{\phi}_T(t)$, we get 
\begin{equation*}
    \Tilde{\phi}_T(t) = \frac{\sqrt{2\pi}T^2}{4t^2T^2 + 4tT^2 + 6\pi t T + T^2 + 3\pi T + 2\pi^2}.
\end{equation*}
After taking the limit as $T \to \infty$, we obtain the expected answer: 
$$\lim_{T \to \infty} \Tilde{\phi}_T(t) = \frac{\sqrt{2\pi}}{4t^2 + 4t + 1} = \frac{\sqrt{2\pi}}{(2t + 1)^2} = h^{\mu}(t).$$

\end{example}

\begin{example}[\cite{PZ}] Consider the PW system with $$d\mu(x) = (1 + \frac{\sin x}{x}) dx.$$ The density of this measure is bounded and bounded away from zero. Therefore $\mu$ is a PW measure. By Theorem \ref{PWConv}, we can use the Hamiltonians recovered from the periodizations $\mu_T$ to approximate the Hamiltonian corresponding to $\mu$. The approximations with $T = \pi, 2\pi, 4\pi, 8\pi$ are in Figure \ref{sincpic}. 

\begin{figure}[ht]
\centering
\includegraphics[width=0.24\textwidth]{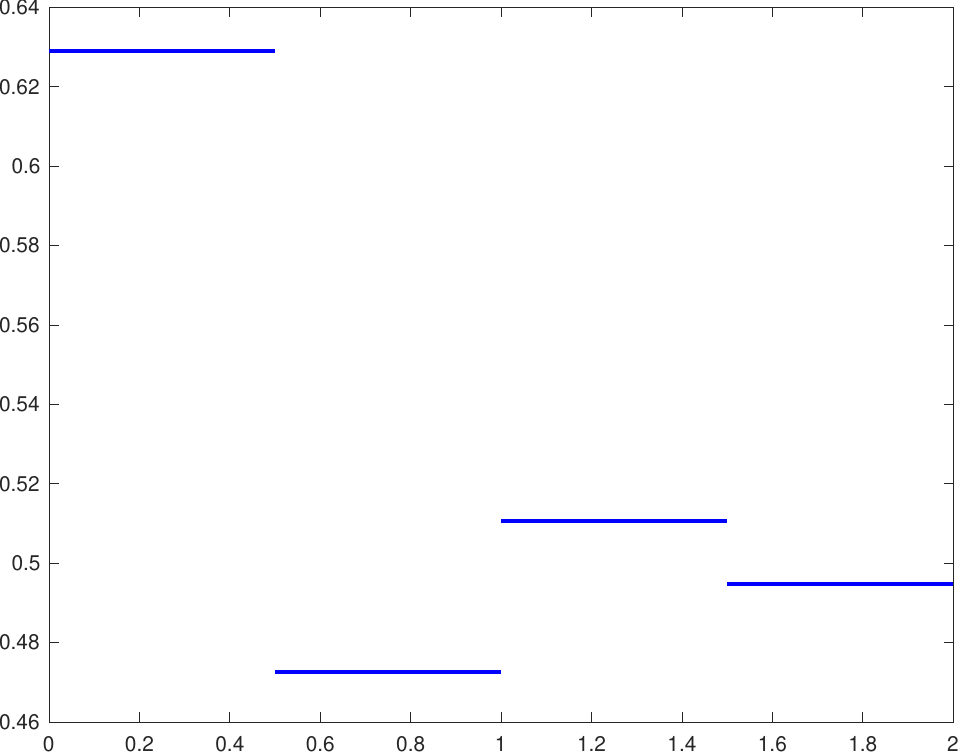}
\includegraphics[width=0.24\textwidth]{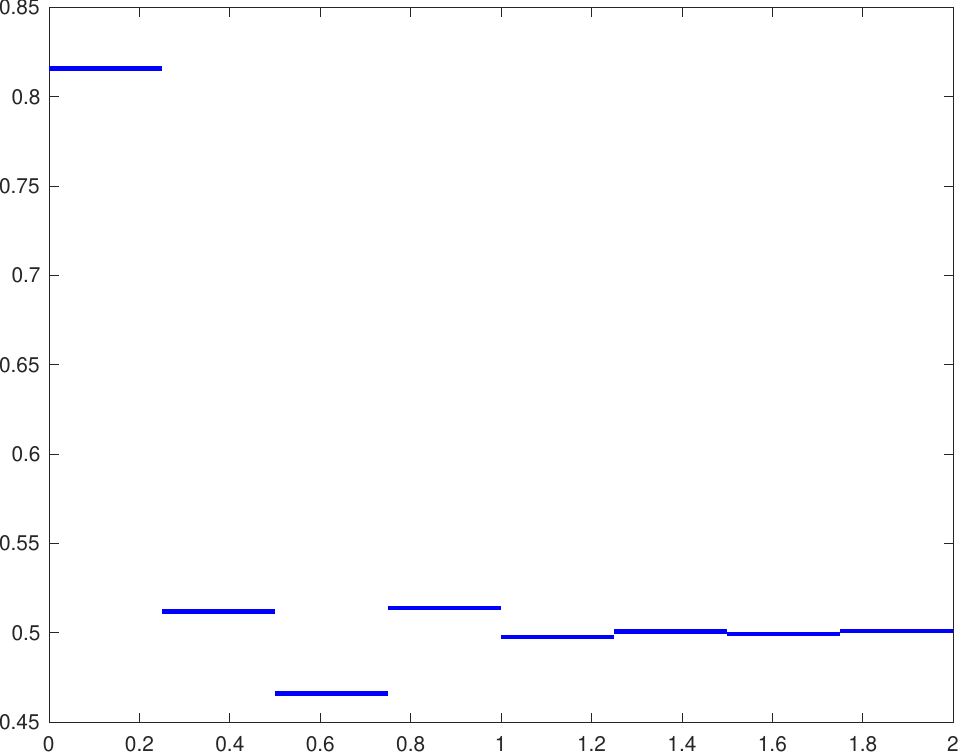}
\includegraphics[width=0.24\textwidth]{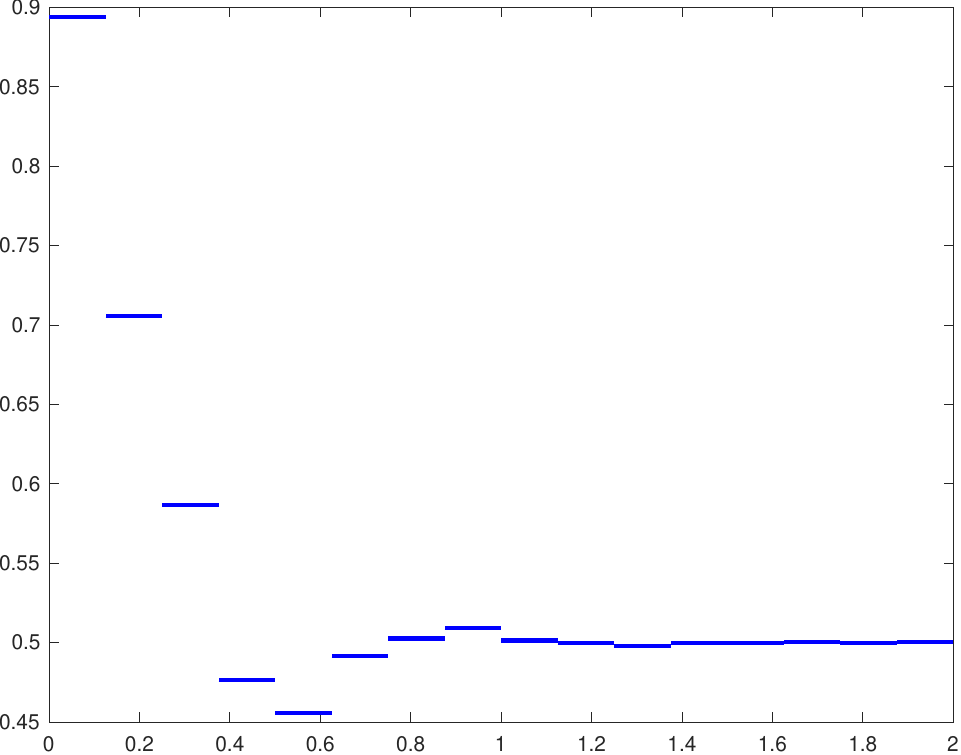}
\includegraphics[width=0.24\textwidth]{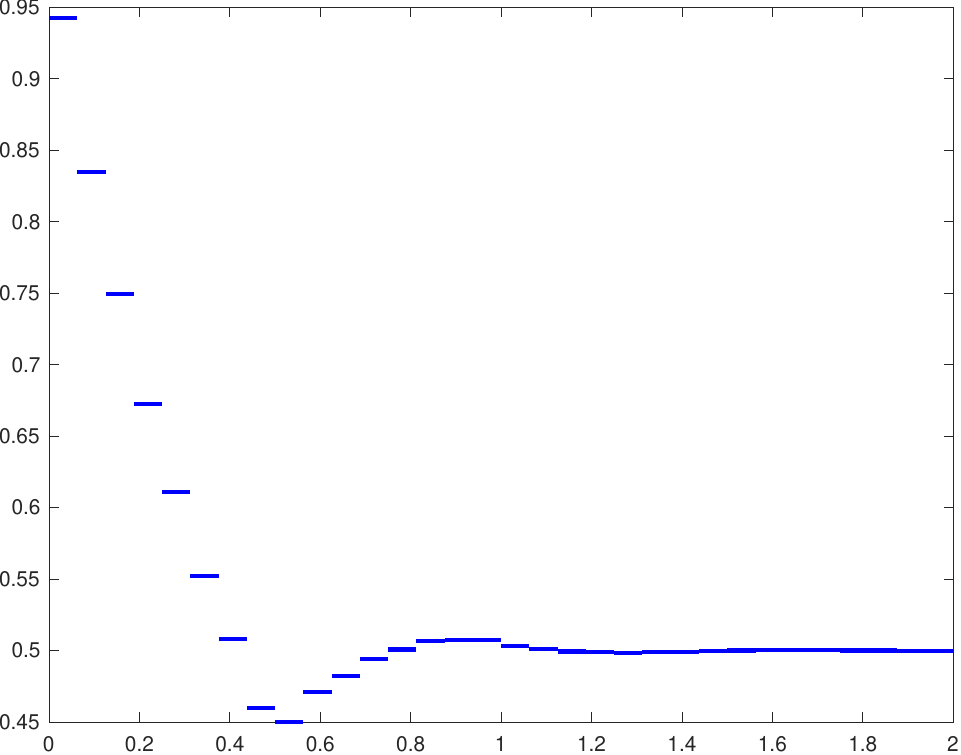}
\caption{Periodization for $d\mu(x) = (1 + \frac{\sin x}{x}) dx$.\\
From left to right: $h^{\mu}_{\pi}$, $h^{\mu}_{2\pi}$, $h^{\mu}_{4\pi}$, $h^{\mu}_{8\pi}$.}
\label{sincpic}
\end{figure}
\end{example}

\begin{example}[\cite{PZ}]
Consider the canonical system with
\[
d\mu(x) = \left( 1 + \sin(x^2) \right) dx.
\]
One can show that $\mu$ is a PW measure using the equivalent conditions in
Section~\ref{PW}. Using Theorem~\ref{PWConv}, we can use the Hamiltonians
recovered from the periodizations $\mu_T$ to approximate the Hamiltonian
corresponding to $\mu$. In particular, we look at $h^{\mu}_T$'s for
$T=\pi, 2\pi, 4\pi, 8\pi$ in Figure~\ref{sinx2pic}.

\begin{figure}[ht]
\centering
\includegraphics[width=0.24\textwidth]{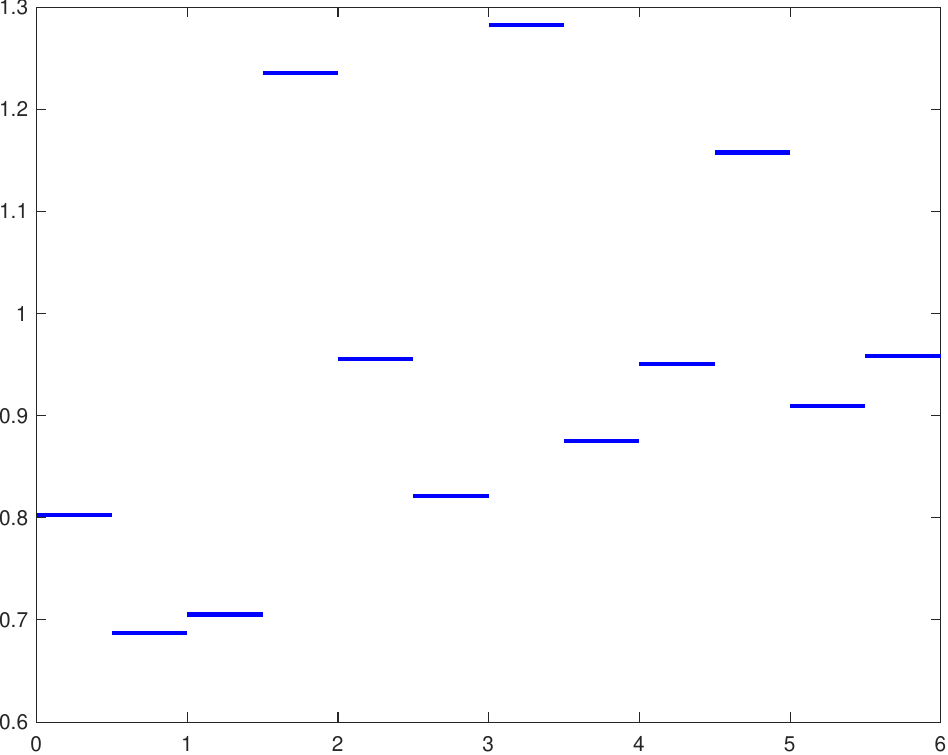}
\includegraphics[width=0.24\textwidth]{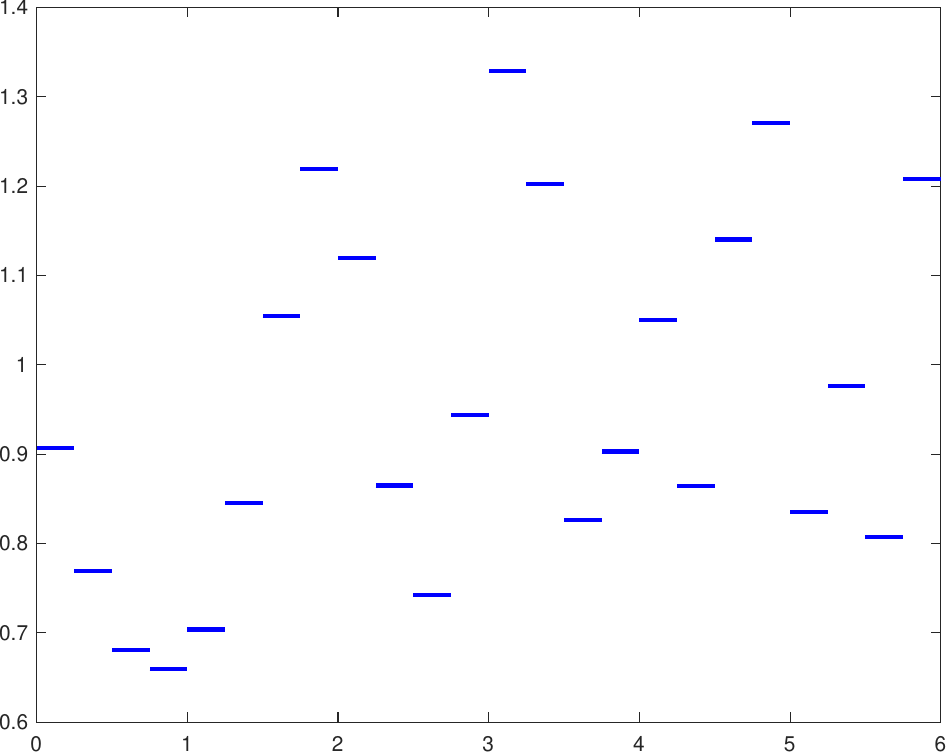}
\includegraphics[width=0.24\textwidth]{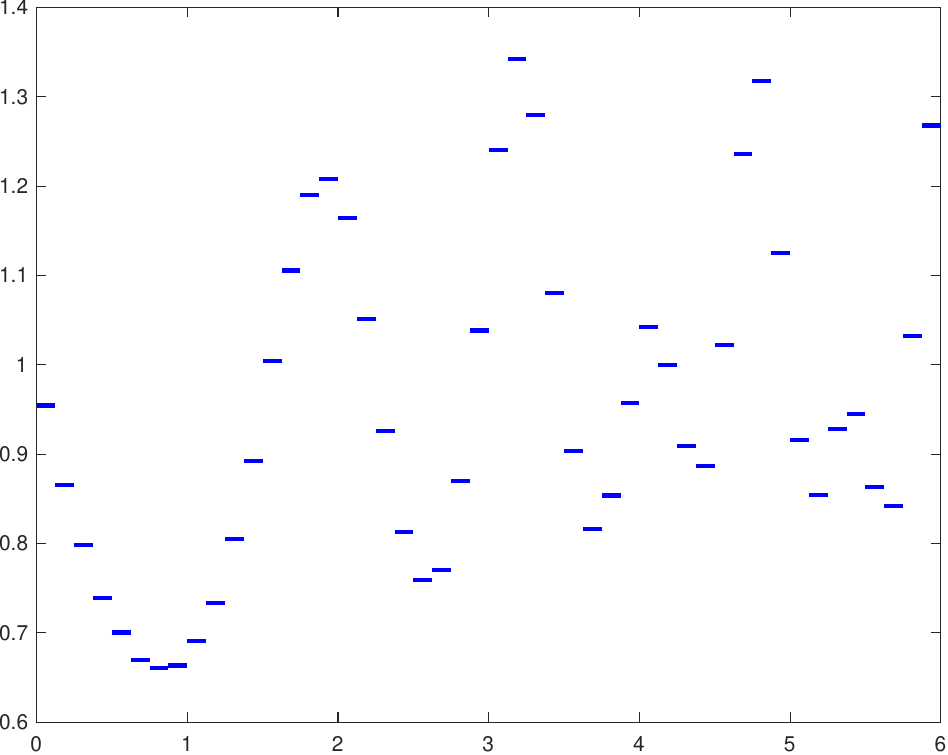}
\includegraphics[width=0.24\textwidth]{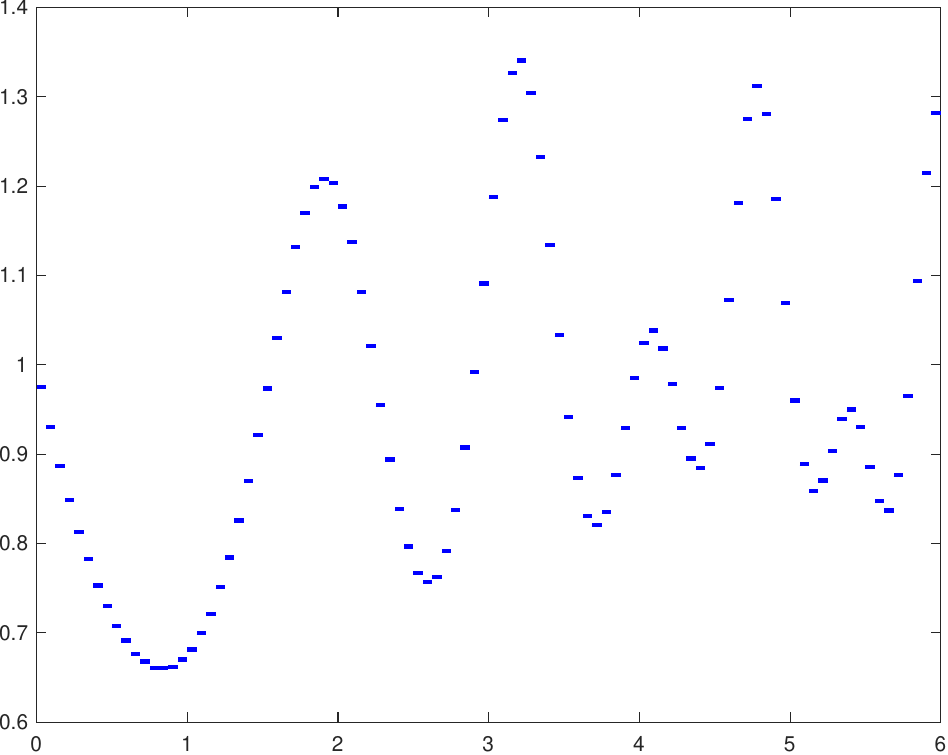}
\caption{Periodization for $d\mu(x) = \left( 1 + \sin(x^2) \right) dx$.\\
From left to right: $h^{\mu}_{\pi}$, $h^{\mu}_{2\pi}$, $h^{\mu}_{4\pi}$, $h^{\mu}_{8\pi}$.}
\label{sinx2pic}
\end{figure}
\end{example}

\begin{example}[\cite{PZ}]
Consider the canonical system with $$d\mu(x) = (1 + |x|^\frac{1}{2})dx.$$ This measure is not a PW measure, but satisfies the assumptions of corollary \ref{CorPolygrowth}.

The approximations with $T = \pi, 2\pi, 4\pi, 8\pi$ can be found in Figure \ref{HalfAbsPic}.

\begin{figure}[ht]
\centering
\includegraphics[width=0.4\textwidth]{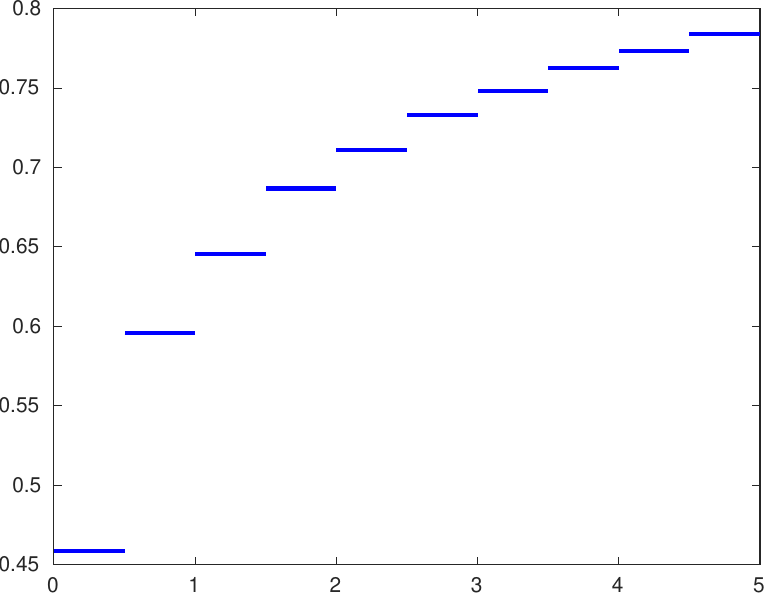}
\includegraphics[width=0.4\textwidth]{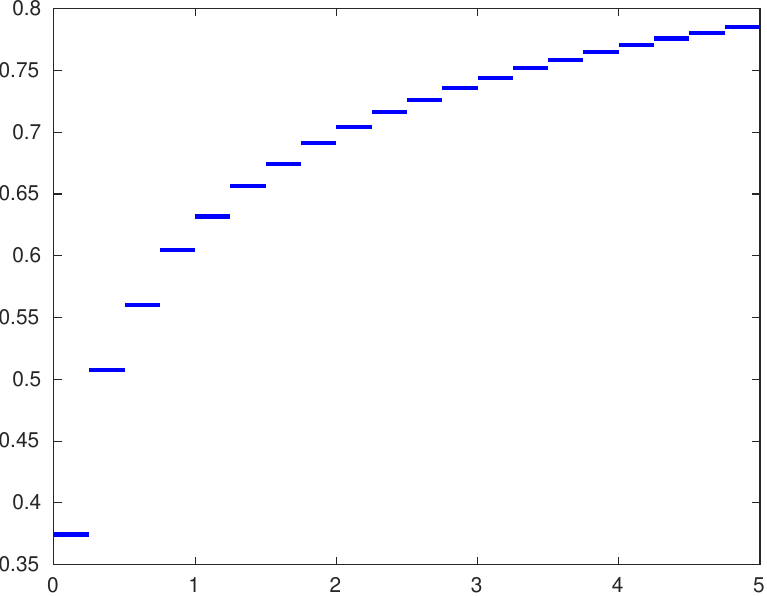}
\includegraphics[width=0.4\textwidth]{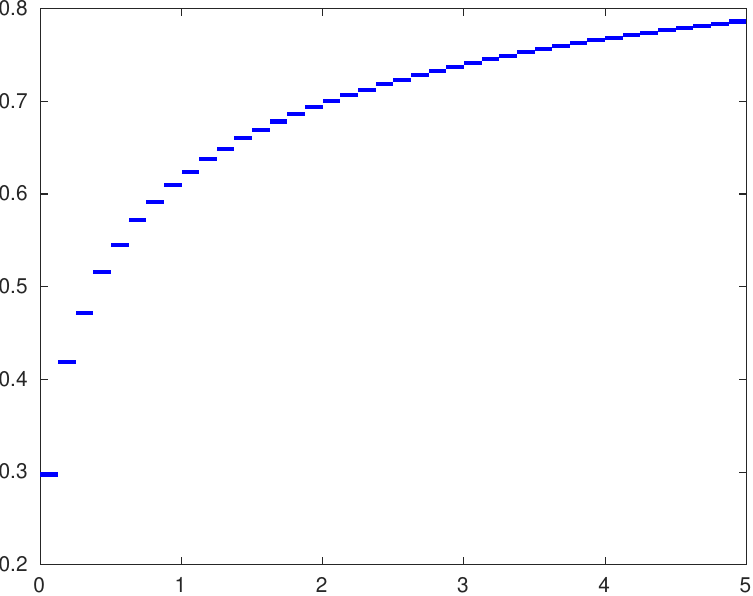}
\includegraphics[width=0.4\textwidth]{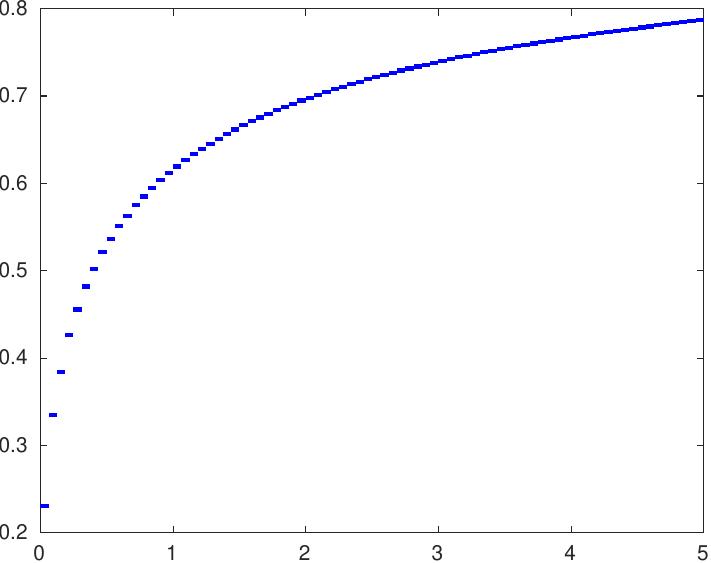}
\caption{Periodization for $d\mu(x) = (1 + |x|^\frac{1}{2})dx$.\\
From left to right: $h^{\mu}_{\pi}$, $h^{\mu}_{2\pi}$, $h^{\mu}_{4\pi}$, $h^{\mu}_{8\pi}$.}
\label{HalfAbsPic}
\end{figure}
\end{example}

\begin{example}[\cite{PZ}]
Consider the canonical system with $$\mu(x) = (1 + |x|)^{\frac{1}{4}}m + \delta_0,$$ where $m$ is the Lebesgue measure and $\delta_0$ is the unit point mass at $0$. This measure is not a PW measure, but satisfies the assumptions of theorem \ref{polygrowth}.

In Figure \ref{DeltaFourthPic}, we present the approximations with $T = \pi, 2\pi, 4\pi, 8\pi$.

\begin{figure}[ht]
\centering
\includegraphics[width=0.4\textwidth]{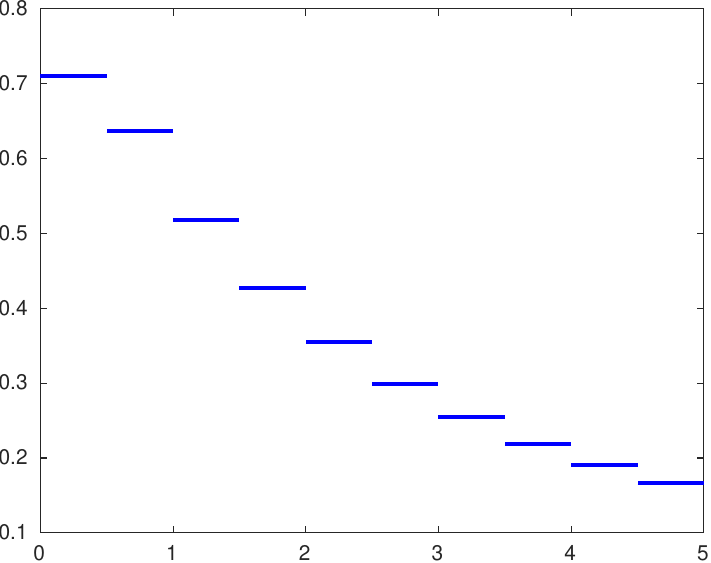}
\includegraphics[width=0.4\textwidth]{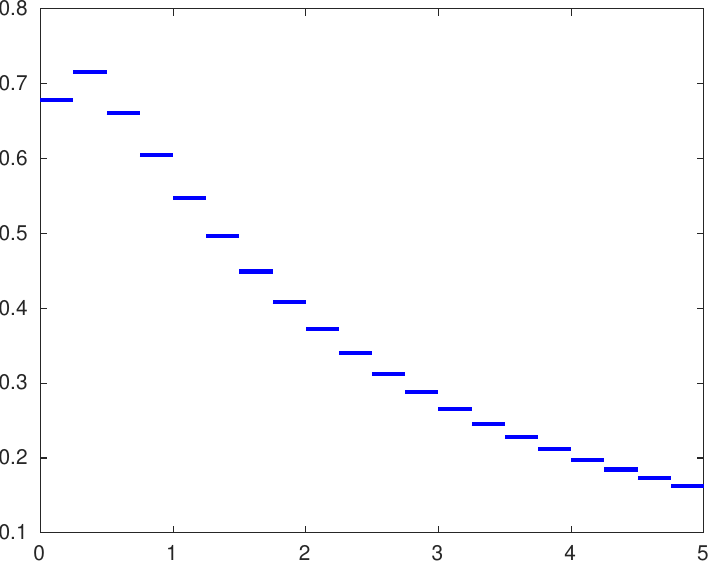}
\includegraphics[width=0.4\textwidth]{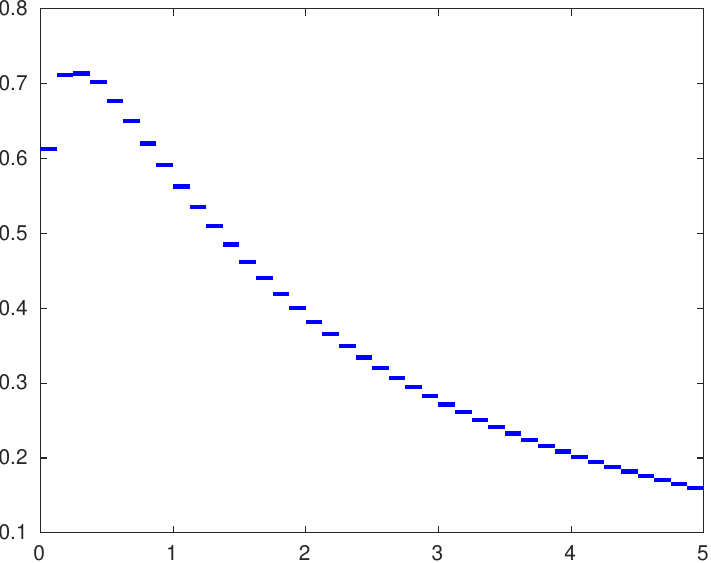}
\includegraphics[width=0.4\textwidth]{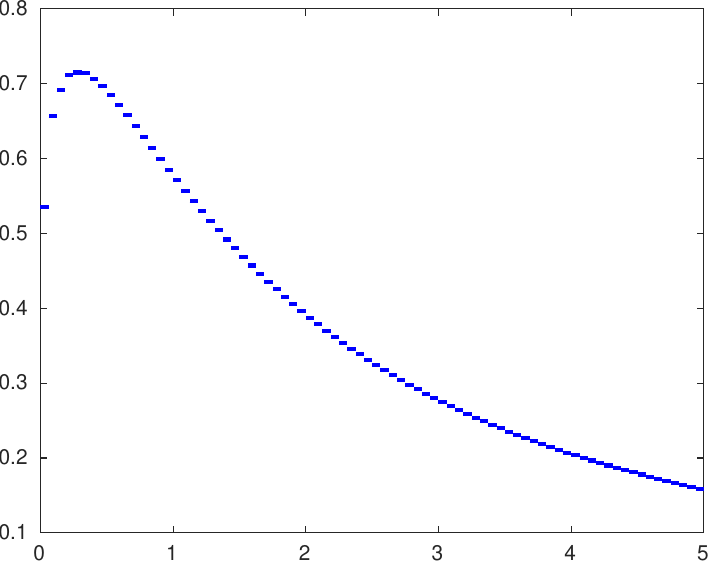}
\caption{Periodization for $\mu(x) = (1 + |x|)^{\frac{1}{4}}m + \delta_0$.\\
From left to right: $h^{\mu}_{\pi}$, $h^{\mu}_{2\pi}$, $h^{\mu}_{4\pi}$, $h^{\mu}_{8\pi}$.}
\label{DeltaFourthPic}
\end{figure}
\end{example}

\begin{example}[\cite{Direct}] \label{ExpGrowth}

Let $H$ be the diagonal det-normalized Hamiltonian with $h_{11}(t) = \exp(t)$. The step-function approximations $h_{11}^T(t)$ satisfy \begin{equation*}
\begin{split}
    \frac{h_{11}^{T, n + 1}}{h_{11}^{T, n}} = \frac{\frac{1}{T} \int_{(n + 1)T}^{(n + 2)T} e^s ds}{\frac{1}{T} \int_{nT}^{(n + 1)T} e^s ds} = e^T,
\end{split}
\end{equation*}
which implies their corresponding spectral measures have constant Verblunsky coefficients.

Using Example \ref{GeometricGrowth}, we can solve for the spectral measures $\mu_T$. For each $T$, on the interval $[0, T\pi] \subseteq \mathbb{R}$, the measure $\mu_T$ has an absolutely continuous part \begin{equation*}
    w^T(\theta) = \frac{T\left( e^T + 1 \right)}{2\left( e^T - 1 \right)} \frac{\sqrt{\frac{4*e^T}{(e^T + 1)^2} - \cos^2(T \theta)}}{\sin(T \theta)},
\end{equation*}
supported on $[\frac{1}{T} \arcsin \frac{1 - e^T}{1 + e^T}, \frac{\pi}{T} - \frac{1}{T} \arcsin \frac{1 - e^T}{1 + e^T}]$, and no singular part.

Since $w^T(x)$ is explicitly known for each $T$, we can determine the limiting measure $\mu$ as well. The limit of $\mu_T$ is given by  
\begin{equation}
    d\mu(x) = \frac{\sqrt{4x^2 - 1}}{2 |x|} dx,
\end{equation}
supported on $(-\infty, -\frac{1}{2}] \cup [\frac{1}{2}, \infty)$.

\end{example}

\subsection{Solution of ISP for a homogeneous system in PW case} \label{ISPHom}

In this section, we present a summary of the calculation from \cite{etudes2} on obtaining solving the ISP for a homogeneous system in the PW case.

Consider homogeneous spectral measures on $\R$ as $$\mu = c_1m + c_2\sigma, $$
with $c_1>|c_2|$, $c_1,c_2\in \mathbb R$. Recall from Section \ref{Hom} that homogeneous measure must be absolutely continuous, $d\mu(x)=\rho(x)dx$,  where
$\rho(x)$  is constant on  $\R_+$ and on $\R_-$.

To compute the functions $h^\mu(t)$ and $g^\mu(t)$ (up to a constant) in the Hamiltonian, recall Theorem \ref{homid}: 
$$k_t(z) = tk_1(tz),$$
in particular,
$$k_t(0) = C_1t,\qquad C_1:= k_1(0).$$
Obviously,
$$h^\mu(t) = \pi \dot k_t(0) = \pi C_1.$$
For the second function we apply the generalized Hilbert transform:
$$g^\mu(t)=\pi \dot l_t(0),\qquad l_t(0)=(H^\mu k_t)(0).$$ 
Recall
\begin{equation*}
\begin{aligned}
    (H^\mu k_t)(0)
    &= \frac{1}{\pi} \int \left[ \frac{k_t(x)-k_t(0)}{x}
        + \frac{x k_t(0)}{1+x^2} \right] \, d\mu(x) \\
    &= A t + \frac{t C_1}{\pi} B,
\end{aligned}
\end{equation*}
where
$$A=\frac 1\pi\int\left[ \frac{k_1(y)-k_1(0)}y+\frac{yk_1(0)}{1+y^2} \right] ~\rho(y)dy$$
and
$$B=\int y\left[\frac1{y^2+t^2}-\frac1{1+y^2}\right]~\rho(y)dy.$$

\noindent Thus $A$ and $C_1$ are constants but $B=B(t)$. In fact,
$$B=-\left[\rho(1)-\rho(-1)\right]\log t$$
and
$$g^\mu(t)=C-\frac 1\pi C_1\left[\rho(1)-\rho(-1)\right]\log t$$
for some constant $C$. 

For det-normalized systems,
$$h_{22}=\frac{1-h_{12}^2}{h^{\mu}}=\frac{1-(C-C_2\log t)^2}{C_1},$$
where $C_2=\frac 1\pi C_1\left[\rho(1)-\rho(-1)\right]$.

This results in the the solution to the ISP, with a free  constant parameter $C$
$$\HH=\begin{pmatrix} C_1& C-C_2\log t\\C-C_2\log t&\frac{1-(C-C_2\log t)^2}{C_1}\end{pmatrix}.$$

As was shown in \cite{etudes}, every spectral measure $\mu$ gives rise to a one-parameter family of Hamiltonians $\HH^\mu$. 
The above formula gives a general solution to the ISP  for $\mu$ defined by \eqref{eq005}.

It remains to compute $C_1$ in terms of $c_1$ and $c_2$. This calculation is technically involved but it highlights its connection to several classical tools of analysis: the Hilbert transform, solutions to the Riemann-Hilbert problem, and special integrals. Here, we present only a brief sketch of the full computation from \cite{etudes2}.

We use the notation $\mathbf{H}f$ for the standard Hilbert transform of a function $f$ on $\R$:
\[
\mathbf{H}f(x) = \frac{1}{\pi}~\text{p.v.}~\int_{\R} \frac{f(t)}{t-x} \, dt.
\]

In terms of convolutions,
\[
\mathbf{H}f = \frac{1}{\pi} f \ast \frac{1}{t}, \quad \text{and} \quad \mathbf{H}\psi_t = \frac{1}{\pi} \, \psi_t \ast \frac{1}{t}.
\]

Let as before $\psi_t=\hat{k_t}$.
Since
$$\psi_t\ast\hat\mu=\sqrt{2\pi}\left[c_1\psi_t-ic_2\mathbf{H}\psi_t\right],$$
from \eqref{eq002A} we obtain
$$\sqrt{2\pi}\left[c_1\psi_t-ic_2\mathbf{H}\psi_t\right]=1\quad {\rm on}\quad(-t,t).$$

To solve this equation for $\psi_t$ we need to deviate into classical complex analysis. We introduce the function
\begin{equation}C_t(z)=\frac1{2\pi i}\int_{-t}^t\frac{\psi_t(s)~ds}{s-z},
\quad z\in \hat{\mathbb C}\setminus [-t,t],\label{eq002}\end{equation}
and denote its boundary values by
$C^\pm_t(x)$, $-t\le x\le t$:
$$C^\pm_t(x)=\lim_{y\to\pm 0}C_t(x+iy).$$ 

\begin{theorem}[\cite{etudes2}]
$$\psi_t=C^+_t-C^-_t,\qquad -i(H\psi_t)=C^+_t+C^-_t.$$
\end{theorem}

It follows that \eqref{eq002A} has the form
\begin{equation}
	C^+_t=GC^-_t+g,\label{eq001}
\end{equation}
where $G$ and $g$ are the numbers
$$G=\frac{c_1-c_2}{c_1+c_2},\qquad g=\frac1{\sqrt{2\pi }}\frac1{c_1+c_2}.$$

The above equation \eqref{eq001} is a particular case of the classical Riemann-Hilbert problem.
Given two functions $F(s)$ and $f(s)$ on $[-t,t]$, we want to find an analytic function $\Phi=\Phi_t(z)$ in $\hat{\mathbb C}\setminus [-t,t]$
such that $\Phi(\infty)=0$ and
$$\Phi^+=F\Phi^-+f\quad {\rm on}\quad [-t,t].$$

The following statement can be verified by direct calculations

\begin{theorem}[\cite{etudes2}] Denote 
\begin{equation}X_t(z)=\exp\left\{\frac1{2\pi i}\int_{-t}^t\frac{\log F(s)}{s-z}~ds\right\}.\label{eq003}\end{equation}
Then
\begin{equation}\Phi_t(z)=X_t(z)\left[\frac1{2\pi i}\int_{-t}^t\frac{f(s)~ds}{(s-z) X_t^+(s)}        \right].\label{eq004}\end{equation}
\end{theorem}

We use the theorem with $F(s)=G$ and $f(s)=g$ to solve \eqref{eq001}.
As $z\to\infty$, \eqref{eq002} implies
$$ C_t(z)~\sim -\frac1z ~\frac1{2\pi i}\int_{-t}^t\psi_t$$
and, from \eqref{eq004},
$$C_t(z)=\Phi_t(z)~\sim -\frac1z ~\frac1{2\pi i}\int_{-t}^t \frac{g}{X_t^+}(s).$$
 
By direct calculations,
$$\int_{-t}^t \frac{1}{s-z}ds=\log\left(\frac{z-t}{z+t}\right).$$
The boundary limits of the last function on $(-t,t)$ are
$$\left(\log\left(\frac{z-t}{z+t}\right)\right)_\pm=\pm\pi i + \log\left|\frac{z-t}{z+t}\right|.$$
From \eqref{eq003} we obtain
$$X^+_t(s)=\exp\left[D\left(\frac 12   +\frac1{2\pi i}\log\left|\frac{s-t}{s+t}\right|\right)\right],$$
where $D=\log G$.

With change of variable $u = \frac{s}{t}$, and further $w = \frac{1 - u}{1 + u}$,  \begin{equation*}
    I=\int_{-t}^t \frac 1{X^+_t(s)}ds = t e^{-D/2} \int_{-1}^{1} \left( \frac{1 - u}{1 + u} \right)^{\frac{D i}{2\pi}} du = t e^{-D/2} \cdot 2 \int_{0}^{\infty} \frac{w^{\frac{D i}{2\pi}}}{(1 + w)^2} dw.
\end{equation*}
Let \( k = \frac{D i}{2\pi} \). Then the last integral becomes
\[
J = \int_{0}^{\infty} \frac{w^{k}}{(1 + w)^2} dw.
\]
This is a Beta integral:
\[
J = B(k + 1, 1 - k) = \frac{\Gamma(k + 1) \Gamma(1 - k)}{\Gamma(2)} = \Gamma(k + 1) \Gamma(1 - k), \quad \text{since } \Gamma(2) = 1.
\]

Hence,
\[
J = \Gamma\left(1 + \frac{D i}{2\pi}\right) \Gamma\left(1 - \frac{D i}{2\pi}\right).
\]

Using the Gamma function identity,
\[
\Gamma(1 + z) \Gamma(1 - z) = \frac{\pi z}{\sin(\pi z)}
\]
with \( z = \frac{D i}{2\pi} \), we obtain
\[
J = \frac{\frac{D i}{2}}{i \sinh\left( \frac{D}{2} \right)} = \frac{D}{2 \sinh\left( \frac{D}{2} \right)}.
\]

Combining the calculations we get
\[
I = t e^{-D/2} \cdot 2 \cdot J = t e^{-D/2} \cdot \frac{D}{\sinh\left( \frac{D}{2} \right)} = \frac{2 t D}{e^{D} - 1}.
\]

To finish the calculation of $C_1$,

\begin{equation*}
    C_1 = \pi \frac{\partial}{\partial t} \int_{-t}^t \psi_t = \pi \frac{\partial}{\partial t}\int_{-t}^{t} \frac{g}{X^+_t(s)} = \frac{2\pi g\log G}{G-1}=\sqrt{\frac\pi 2}\frac 1{c_2}\log\left(\frac{c_1+c_2}{c_1-c_2}\right).
\end{equation*}
 
And substituting this into $C_2$ we get
$$
C_2= \frac{1}{\sqrt{2\pi}}\log\left(\frac{c_1+c_2}{c_1-c_2}\right).
$$

\subsection{Bessel functions and canonical systems}\label{bessel}
In this section we present the example of a direct spectral problem involving families of Bessel functions from \cite{etudes2}.

The Bessel functions can be defined using the following lemma
\begin{lemma}[\cite{etudes2}] Suppose $F(t)$ satisfies 
\begin{equation}t^2\ddot F(t)+t\dot F(t)+(t^2-\nu^2)F(t)=0\label{eq007}\end{equation}
for $t>0$. Let $\kappa>0, \beta>0$, $\alpha$ be given real numbers. Then the function
$$y(t)=t^\alpha F(\kappa t^\beta) $$
solves the equation
\begin{equation}\qquad t^2 \ddot y+at\dot y+(b+c^2t^{2\beta})y=0\label{eq006}\end{equation}
with
\begin{equation}\qquad a=1-2\alpha,\qquad  b=\alpha^2-\beta^2\nu^2, \qquad c^2=\beta^2\kappa^2.\label{eq008}\end{equation}
\end{lemma}

Let $J_\nu$ denote the Bessel function of the first kind. Then $J_\nu$ solves \eqref{eq007} for all $t\in \R$ and is finite at $0$, which can be taken as the definition of $J_\nu$.

The special case that we study in this section is
$$t^2\ddot y+at\dot y+c^2t^2 y=0$$
with parameters $a$ and $c>0$ (and also  $\beta=1$, $b=0$).  We find 
$$\alpha=\frac{1-a}2,\qquad \nu=\alpha,\qquad \kappa=c.$$ 
The general solutions is
$${\rm span} \left\{t^\alpha J_\alpha(c t),\;
t^\alpha J_{-\alpha}(c t)\right\}.$$

We are now ready to define the Bessel canonical system. For each $m>0$ we consider the system with the Hamiltonian
$$H(t)=\begin{pmatrix}
   h(t)   & 0 \\
     0& 1/ h(t)
\end{pmatrix}, \qquad h(t)=t^m,\quad t>0.$$
Note that the system is {\it regular} iff $m<1$.  As usual we introduce the functions $A=A(t,z)$ and $C=C(t,z)$:
$$\Omega\dot X=zHX,\qquad X=(A, C)^\tau,$$
i.e.,
$$\dot C=zhA,\quad -\dot A=\frac z h C,$$
with "initial values"
$$A(0, z)=1, \qquad C(0, z)=0$$
(as limits when $t\to 0$). 
Rewriting the system as a second order equation for $C$ we get
$$\ddot C-\frac mt\dot C+z^2C=0$$ 
with initial conditions
$$C(0,z)=0,\qquad \dot C(t,z)\sim zt^m\text{ as }t\to 0.$$

Consider the function $F_\nu$ defined by
$$J_\nu(\lambda)=\lambda^\nu F_\nu(\lambda).$$
It is known that $F_\nu$ is an entire function and
$$F_\nu(0)=\frac1{2^\nu\Gamma(\nu+1)},\qquad F'_\nu(0)=0.$$
If we fix $z$ and let $t\to 0$, then we have
$$t^\nu J_\nu(zt)=t^{2\nu}z^\nu F_\nu(zt)\sim F_\nu(0) z^\nu t^{2\nu}$$
and 
$$t^\nu J_{-\nu}(zt)=z^{-\nu }F_{-\nu}(zt) \to z^{-\nu}F_{-\nu}(0).$$
From the condition $C(0,z)=0$ it follows that
$$C(t,z)=G(z)t^\nu J_\nu(zt),$$
It remains to find $G(z)$. 

We have 
$$C(t,z) \sim G(z)  F_\nu(0) z^\nu t^{2\nu},$$
and
$$ \dot C(t,z) \sim G(z)  F_\nu(0) z^\nu ~ 2\nu t^{2\nu-1}.$$
Combining with the second boundary condition, 
$
\dot C(t,z)\sim zt^m$, we have
$$zt^m\sim G(z)  F_\nu(0) z^\nu ~ 2\nu t^{2\nu-1},$$
and therefore (recall that $m=2\nu-1$)
$$ G(z)=\frac{z^{1-\nu}}{F_\nu(0) ~2\nu}=g_\nu z^{1-\nu},\qquad g_\nu:=2^{\nu-1}\frac{\Gamma(1+\nu)}\nu=2^{\nu-1}\Gamma(\nu).$$

Thus
$$C(t,z)=g_\nu t^{2\nu}z F_\nu(zt).$$
For each fixed $t$, $C$ is entire with respect to $z$. We arrive at
$$C=g_\nu(t^{2\nu} z F_\nu(zt)),$$
and 
$$A=g_\nu(2\nu F_\nu(zt)+tz F'_\nu(zt)).$$
Indeed, we have
$$\dot C=g_\nu(2\nu t^{2\nu-1} zF_\nu(zt)+t^{2\nu}z^2 F'_\nu(zt)),$$
and, since $m=2\nu-1$,
$$A=g_\nu\frac{\dot C}{zt^m}=g_\nu(2\nu F_\nu(zt)+tz F'_\nu(zt)).$$ 

We can simplify our expression for $A$ to obtain:

\begin{theorem}[\cite{etudes2}]
$$A(t,z)=g_\nu F_{\nu-1}(zt)$$
and
$$C(t,z)=g_\nu t^{2\nu}z F_\nu(zt).$$
\end{theorem}
 
Let $k_t(z)$ be the reproducing kernel of $\BB(E_t)$, $E:=A-iC$, at zero:
 $$k_t(z)=\frac{C(z)A(0)}{\pi z}=\frac{g_\nu^2 F_{\nu-1}(0)}\pi t^{2\nu} F_\nu(zt).$$
 
Together with the resulsts of  \cite{etudes2}, this relation implies that the spectral measure $\mu$ of the system is quasi-homogeneous of order $\nu-1=(m-1)/2$. The diagonal form of the Hamiltonian implies that $\mu$ is even. Altogether, 
 $$\mu=\const |x|^{m}.$$
 
As we can see, $\mu\not\in \PW$.

In conclusion, let us verify our computations by 'deriving' the Hamiltonian. We can apply \eqref{eq003A}, 
even though the system does not satisfy the PW requirement.
The correctness of the answer obtained suggests broader use of Theorem \ref{t3} and similar formulas.
 
According to our previous calculations,
 $$\pi k_t(0)=\frac{g_\nu^2 F_{\nu-1}(0)}\pi t^{2\nu} F_\nu(0)= \frac{1}{2\nu\pi} t^{2\nu}.$$ 
 Via \eqref{eq003A},
 $$h^{\mu}(t)=\pi\dot k_t(0)=t^{2\nu-1}=t^m.$$

\subsection{Problem: Description of defining sequences in the case of condition-free endpoint in the two-interval problem.}\label{P1}

In Horvath' theorem discussed in Section \ref{sHorvath} the operator is recovered from a part of its potential on $(0,a\pi)$ and
the values of the Weyl function $m_-$, obtained by fixing a boundary condition at the opposite point $\pi$. What if
we consider a similar problem with a condition-free end at $\pi$, i.e., try to recover an operator from its potential near $0$ and
the function $m_+$ obtained by fixing a boundary condition at the same endpoint $0$? In this case the statement
analogous to Horvath' theorem fails.

\begin{counterexample}
The simplest counterexample is a modification of the 3-interval counterexample from \cite{SG2} mentioned before.
Let the operators $L$ and $\ti L$  from $\SS^2$ be such that
$q=\ti q \textrm{ on }(0,(1-\e)\pi)$ and $q(x)=\ti q((2-\e)\pi-x)$ on $((1-\e)\pi,\pi)$. Consider an auxiliary operator
$L_a$ on $(0,(2-\e)\pi)$ whose potential is defined as $q_a=q$ on $(0,\pi)$ and $q_a(x)=q((2-\e)\pi)-x)$ on $(\pi,(2-\e)\pi)$.
Let $\ti L_a$ be an operator such that $\ti q_a=\ti q$ on $(0,\pi)$ and $\ti q_a(x)=q((2-\e)\pi-x)$ on $(\pi,(2-\e)\pi)$.
Note that then the DD spectra of $L_a$ and $\ti L_a$ coinside.
Denote this sequence by $\L$. It is not difficult to see the Weyl functions $m_+$ and $\ti m_+$ of the operators $L$ and $\ti L$
will coincide on $\L$. It is left to notice that the density of $\L$ is $2-\e$, i.e., much larger than $\e$.
\end{counterexample}

Although, as we can see from this counterexample, not any sequence $\L$ of large enough density can be used in our modified two-interval problem, there are
some sequences which can  be used, as shown by the following.

\begin{proposition}[\cite{MP3}]
Let $g\in L^2((0,\pi-\e)),\ \e<\pi/2$ and let $\Sigma=\sigma_{DD}$ be the spectrum of  the \Sch operator from $\SS^2((0,\pi-\e))$ with $q=g$.
Then any two $L,\ti L\in \SS^2((0,\pi))$, such that $q=\ti q=g$ on $(0,\pi-\e)$ and $m_+=\ti m_+$ on
a subsequence $\L$ of $\Sigma$, $\pi D^*(\L)>2\e$, must be identical.
\end{proposition}

Let $g\in L^2(0,\pi-\e)$. Denote by $\SS^2_g$ the set of all operators from $\SS^2$ such that $q=g$ on $(0,\pi-\e)$.
The natural problem which arises from the proposition is to describe the set of sequences $\L$ such that
the values of the Weyl function $m_+$ on $\L$ determine an operator uniquely within the class $\SS^2_g$. 
The counterexample shows that the set does not contain all of the sequences of proper density and the Proposition
shows that the set is non-empty.

\subsection{Problem: Description of unique operators in the three-interval case.}\label{P2}

Let $0\leq a<b\leq 1$ and $a+(1-b)>1/2$.
As follows from the discussion in Section \ref{3intU}, there exist \Sch operators  $L\in \SS^2$ on $(0,\pi)$ with the following
uniqueness property: If $\ti L\in \SS^2$ is another \Sch operator such that $\ti q=q$ on $(0,a\pi)\cup (b\pi,\pi)$
and $\pi D^*(\ti\sigma_{DD}\cap \sigma_{DD})>2(b-a)$ then $\ti L \equiv L$.

This rises a natural question of description of all such 'unique' operators, i.e., operators
uniquely determined by the restriction of their potentials on $(0,a\pi)\cup (b\pi,\pi)$ and
a properly sized subsequence of the spectrum. The example provided by Theorem \ref{tUoper} presents an operator
which is close to even in the spectral case. The question one may start with is whether there are other examples.

While our discussion in \cite{MP3} concentrates on the $\SS^2$ case, it can be extended to regular case without much difficulty. Further cases
of this problem may concern similar examples and descriptions in non-regular cases. If $q$ is unsummable, the correct question would be to
describe $L$ such that any $\ti L$ with the above properties, and such that $q-\ti q$ is small (summable, for instance), must coincide with $L$.

The three-interval case is only a model case in which we already see some of the difficulties that were not present in the two-interval problem.
The ultimate goal in this and similar problems would be to consider other subsets of the interval.

\subsection{Problem: Uncertainty in other types of spectral data}\label{P3a}

The version of the two-spectra  problem treated in Section \ref{2int} is not the only case when a problem of Uncertainty Quantification appears naturally in spectral settings.  Let us give another example of such a problem.

Let $L\in \SS^2$ be a \Sch operator with the spectral measure $\mu_-=\sum \alpha_n\delta_{\lan}$. Let
$\EE=\{\e_n\}$ be a sequence of positive numbers. Denote by $\EE_L$ the set of operators $\ti L\in \SS^2$
such that $\ti\mu_- =\sum\ti \alpha_n\delta_{\lan}$, $|\alpha_n-\ti \alpha_n|\leq\e_n$. Similarly
to Section \ref{2int}, denote by $U(\EE_L)$ the infimum of $a$ such that the values of the potential on $(0,a\pi)$
determine a \Sch operator from $\EE_L$ uniquely.

\begin{proposition}[\cite{MP3}]
$$U(\EE_L)=\pi \sup \left\{ D_*(\Phi)\ :\ \sum_{n\in\Phi}\frac{\log_-\e_n}{1+n^2}<\infty\right\}.$$
\end{proposition}

The proof follows easily from Lemma \ref{l5} and Corollary \ref{l2}.

The estimates of the size of uncertainty in other variations of spectral problems may require different techniques.
The definition of $U(D)$ for the spectral data $D$ which worked for us in the problems considered in this paper may need to be further developed for other kinds of data. Let us point out, for instance, that the present
definition of $U$  will not work properly if one replaces $\mu_-$ with $\mu_+$ in the definition of $\EE_L$ above. Even though
an obvious modification of $U$ by replacing  $(0,a\pi)$ with $((1-a)\pi,\pi)$ will fix this particular problem,
one would like to have a more universal definition.

\subsection{Problem: Further connections between mixed spectral problems and completeness problems.}\label{P3}

Horvath' Theorem \ref{Horvath} has established a connection between the Beurling-Malliavin (BM) problem on completeness of systems of exponentials
in $L^2$ on an interval and the original case of the mixed spectral problem, the two-interval case without a condition-free endpoint.
Even though cases of multiple intervals were considered in the literature (see for instance
\cite{Tru}, Chapter 4, Problem 10a) similar connections are yet to be found.
It is interesting to observe that the analog of the BM theorem with one interval replaced with any other subset of the line, including
the next simplest case of a union of two intervals, does not exist. Despite a number of deep results on completeness of exponential systems in $L^2$ over general sets, see for
instance \cite{Olevsky}, there is no formula for the radius of completeness or  a good idea on what could replace such a formula, even in the case of two intervals.

As we saw in Section \ref{s3int}, many of the same complications appear in mixed spectral problems when moving from the two-interval to the three-interval and multiple-interval case.
Moreover, analogies between the Weyl transform and the Fourier transform, together with the use of the latter in BM theory,
suggest that the mixed spectral problems for more general subsets of the interval must be closely related to BM problems over general sets.
Finding such connections, i.e., formulating an analog of Horvath' theorem for more general subsets seems to be an interesting and challenging problem.


\newpage
\begin{thebibliography}{99}


\bibitem{BBP}
{\sc Baranov, A., Belov, Yu., Poltoratski, A.}
{\it De Branges functions for Schr\"odinger equations,}
Preprint, 2016.

\bibitem{BM1}
{\sc Beurling, A., Malliavin, P.}
{\it On Fourier transforms of measures with compact support,}
Acta Math. 107 (1962), 291--302.

\bibitem{BM2}
{\sc Beurling, A., Malliavin, P.}
{\it On the closure of characters and the zeros of entire functions,}
Acta Math. 118 (1967), 79--93.

\bibitem{Bessonov}
{\sc Bessonov, R.}
{\it Sampling measures, Muckenhoupt Hamiltonians, and triangular factorization,}
Int. Math. Res. Not. 2018, no.~12, 3744--3768.

\bibitem{BR}
{\sc Bessonov, R., Romanov, R.}
{\it An inverse problem for weighted Paley--Wiener spaces,}
Inverse Problems 32 (2016), no.~11.

\bibitem{Borg}
{\sc Borg, G.}
{\it Uniqueness theorems in the spectral theory of $y''+(\lambda-q(x))y=0$,}
Proc. 11th Scandinavian Congress of Mathematicians,
Johan Grundt Tanums Forlag, Oslo, 1952, 276--287.

\bibitem{dB}
{\sc De Branges, L.}
{\it Hilbert spaces of entire functions,}
Prentice--Hall, Englewood Cliffs, NJ, 1968.

\bibitem{DGS}
{\sc del Rio, R., Gesztesy, F., Simon, B.}
{\it Inverse spectral analysis with partial information on the potential. III.
Updating boundary conditions,}
Internat. Math. Res. Notices 1997, no.~15, 751--758.

\bibitem{SG3}
{\sc Gesztesy, F., Simon, B.}
{\it Uniqueness theorems in inverse spectral theory for one-dimensional
Schr\"odinger operators,}
Trans. Amer. Math. Soc. 348 (1996), no.~1, 349--373.

\bibitem{SG1}
{\sc Gesztesy, F., Simon, B.}
{\it A new approach to inverse spectral theory, II.
General real potentials and the connection to the spectral measure,}
Ann. of Math. 152 (2000), 593--643.

\bibitem{SG2}
{\sc Gesztesy, F., Simon, B.}
{\it Inverse spectral analysis with partial information on the potential, II.
The case of discrete spectrum,}
Trans. Amer. Math. Soc. 352 (2000), 2765--2787.

\bibitem{HL}
{\sc Hochstadt, H., Lieberman, B.}
{\it An inverse Sturm--Liouville problem with mixed given data,}
SIAM J. Appl. Math. 34 (1978), 676--680.

\bibitem{Horvath}
{\sc Horv\'ath, M.}
{\it Inverse spectral problems and closed exponential systems,}
Ann. of Math. 162 (2005), no.~2, 885--918.

\bibitem{LS}
{\sc Levitan, B.~M., Sargsjan, I.~S.}
{\it Sturm--Liouville and Dirac operators,}
Kluwer, Dordrecht, 1991.

\bibitem{Lub}
{\sc Lubinsky, D.~S.}
{\it A Survey of Weighted Polynomial Approximation with Exponential Weights,}
Surveys in Approximation Theory 3 (2007), 1105.

\bibitem{MIF1}
{\sc Makarov, N., Poltoratski, A.}
{\it Meromorphic inner functions, Toeplitz kernels, and the uncertainty principle,}
in {\it Perspectives in Analysis},
Springer, Berlin, 2005, 185--252.

\bibitem{MP3}
{\sc Makarov, N., Poltoratski, A.}
{\it Two spectra theorem with uncertainty,}
J. Spectral Theory 9 (2019), no.~4, 1249--1285.

\bibitem{etudes}
{\sc Makarov, N., Poltoratski, A.}
{\it \'Etudes in the inverse spectral problem,}
J. London Math. Soc. 108 (2023), no.~3, 916--977.

\bibitem{etudes2}
{\sc Makarov, N., Poltoratski, A.}
{\it \'Etudes in the inverse spectral problem, II,}
Preprint, arXiv, 2025.

\bibitem{MPS}
{\sc N. Makarov, A. Poltoratski, and M. Sodin.}
{"Lectures on linear complex analysis",}
{\it unpublished lecture notes}.

\bibitem{M1}
{\sc Marchenko, V.}
{\it Some questions in the theory of one-dimensional linear differential
operators of the second order, I,}
Trudy Mosk. Mat. Obshch. 1 (1952), 327--420.

\bibitem{M2}
{\sc Marchenko, V.}
{\it Sturm--Liouville operators and applications,}
Birkh\"auser, Basel, 1986.

\bibitem{Meg}
{\sc Mergelyan, S.}
{\it Weighted approximation by polynomials,}
Uspekhi mat. nauk 11 (1956), 107--152,
English translation in Amer. Math. Soc. Translations, Ser. 2, 10 (1958), 59--106.

\bibitem{Polya}
{\sc Mitkovski, M., Poltoratski, A.}
{\it Polya sequences, Toeplitz kernels and gap theorems,}
Adv. Math. 224 (2010), 1057--1070.

\bibitem{Det}
{\sc Mitkovski, M., Poltoratski, A.}
{\it On the determinacy problem for measures,}
Invent. Math. 202 (2015), 1241--1267.

\bibitem{OS}
{\sc Ortega-Cedr\'a, J., Seip, K.}
{\it Fourier frames,}
Ann. of Math. 155 (2002), 789--806.

\bibitem{Olevsky}
{\sc Olevskii, A., Ulanovskii, A.}
{\it Functions with disconnected spectrum: sampling, interpolation, translates,}
Univ. Lecture Series, AMS, 2016.

\bibitem{Gap}
{\sc Poltoratski, A.}
{\it Spectral gaps for sets and measures,}
Acta Math. 208 (2012), no.~1, 151--209.

\bibitem{Type}
{\sc Poltoratski, A.}
{\it A problem on completeness of exponentials,}
Ann. of Math. 178 (2013), 983--1016.

\bibitem{CBMS}
{\sc Poltoratski, A.}
{\it Toeplitz Approach to Problems of the Uncertainty Principle,}
CBMS Series, Amer. Math. Soc., Providence, RI, 2015.

\bibitem{PS}
{\sc Poltoratski, A., Sarason, D.}
{\it Aleksandrov--Clark measures,}
Recent advances in operator-related function theory,
Contemp. Math. 393, Amer. Math. Soc., Providence, RI, 2006, 1--14.

\bibitem{PZ}
{\sc Poltoratski, A., Zhang, A.~R.}
{\it Periodic approximations in inverse spectral problems for canonical
Hamiltonian systems,}
J. Funct. Anal. 284 (2023), no.~11.

\bibitem{Tru}
{\sc P\"oschel, J., Trubowitz, E.}
{\it Inverse spectral theory,}
Academic Press, New York, 1987.

\bibitem{Rem}
{\sc Remling, C.}
{\it Spectral theory of canonical systems,}
De Gruyter, 2018.


\bibitem{Rom}
{\sc R. Romanov.}
{"Canonical Systems and de Branges Spaces",}
2014. \url{arXiv:1408.6022}.

\bibitem{S}
{\sc Simon, B.}
{\it A new approach to inverse spectral theory, I.
Fundamental formalism,}
Ann. of Math. 150 (1999), 1029--1057.

\bibitem{Simon}
{\sc Simon, B.}
{\it Orthogonal polynomials on the unit circle, Parts I and II,}
AMS Colloquium Publications, Vol.~54, 2004.


\bibitem{W}
{\sc Winkler, H.}
{\it On Transformation of Canonical Systems,}
in Gohberg, I., Langer, H. (eds),
{\it Operator Theory and Boundary Eigenvalue Problems},
Operator Theory: Advances and Applications, vol.~80,
Birkh\"auser, Basel.

\bibitem{Zhang}
{\sc Zhang, A.}
{"Convergence from the discrete to the continuous non-linear Fourier transform",}
{\it Transactions of the American Mathematical Society} 378 (2025), 1483–1501.

\bibitem{Direct}
{\sc Zhang, A.}
{\it Direct spectral problems for Paley-Wiener canonical systems,}
\url{https://arxiv.org/abs/2505.00669}.

\bibitem{h12}
{\sc Zhang, A.}
{\it Direct spectral problems for Paley-Wiener canonical systems, II,}
Preprint, 2026.


 

\end{thebibliography}
\end{document}